\documentclass[12pt]{amsart}
\usepackage[all]{xy}

\usepackage{graphicx}  

\usepackage{amsmath}
\usepackage{amssymb}
\usepackage{amsfonts}
\usepackage{mathrsfs}
\usepackage{mathtools}

\newcommand{\Cdb}{\ensuremath{\mathbb{C}}}
\newcommand{\Ddb}{\ensuremath{\mathbb{D}}}

\newcommand{\Ndb}{\ensuremath{\mathbb{N}}}

\newcommand{\Rdb}{\ensuremath{\mathbb{R}}}

\newcommand{\Tdb}{\ensuremath{\mathbb{T}}}

\newcommand{\Zdb}{\ensuremath{\mathbb{Z}}}

\renewcommand{\H}{\mbox{${\mathcal H}$}}

\renewcommand{\O}{\mbox{${\mathcal O}$}}
\renewcommand{\P}{\mbox{${\mathcal P}$}}

\newcommand{\R}{\mbox{${\mathcal R}$}}

\newcommand{\norm}[1]{\Vert#1\Vert}
\newcommand{\bignorm}[1]{\bigl\Vert#1\bigr\Vert}
\newcommand{\Bignorm}[1]{\Bigl\Vert#1\Bigr\Vert}
\newcommand{\biggnorm}[1]{\biggl\Vert#1\biggl\Vert}

\newcommand{\Rad}{{\rm Rad}}

\newtheorem{theorem}{Theorem}[section]
\newtheorem{lemma}[theorem]{Lemma}
\newtheorem{corollary}[theorem]{Corollary}
\newtheorem{proposition}[theorem]{Proposition}
\newtheorem{definition}[theorem]{Definition}
\theoremstyle{remark}
\newtheorem{remark}[theorem]{\bf Remark}
\theoremstyle{definition}

\numberwithin{equation}{section}

\newcommand{\Gauss}{{\rm Gauss}}
\newcommand{\bnorm}[1]{ \big\| #1  \big\|}

\newcommand{\bgnorm}[1]{ \bigg\| #1  \bigg\|}

\newcommand{\ot}{\otimes}

\def\Diag{{\rm Diag} \, }

\def\var{{\rm var} \, }

\begin{document}

\title[]{Dilation of Ritt operators on $L^p$-spaces}

\author{C\'edric Arhancet, Christian Le Merdy}
\address{Laboratoire de Math\'ematiques\\ Universit\'e de  Franche-Comt\'e
\\ 25030 Besan\c con Cedex\\ France}
\email{clemerdy@univ-fcomte.fr}

\address{Laboratoire de Math\'ematiques\\ Universit\'e de  Franche-Comt\'e
\\ 25030 Besan\c con Cedex\\ France}
\email{cedric.arhancet@univ-fcomte.fr}

\date{\today}

\begin{abstract}
For any Ritt operator $T\colon L^p(\Omega)\to L^p(\Omega)$, 
for any positive real number $\alpha$,
and for any $x\in L^p(\Omega)$, 
we consider $\norm{x}_{T,\alpha}\,=\, \Bignorm{\Bigl(\sum_{k=1}^{\infty} k^{2\alpha -1}\bigl 
\vert T^{k-1}(I-T)^\alpha x\bigr\vert^2\Bigr)^{\frac{1}{2}}}_{L^p}\,$.
We show that if $T$ is actually an $R$-Ritt operator, then the square functions
$\norm{\ }_{T,\alpha}$ are pairwise equivalent. Then we show that $T$ and
its adjoint $T^*\colon L^{p'}(\Omega)\to L^{p'}(\Omega)$ both  
satisfy uniform estimates $\norm{x}_{T,1}
\lesssim \norm{x}_{L^p}$ and $\norm{y}_{T^*,1}
\lesssim \norm{y}_{L^{p'}}$ for $x\in L^p(\Omega)$ and $y\in L^{p'}(\Omega)$ if and only if $T$
is $R$-Ritt and admits a dilation in the following sense: there exist a measure
space $\widetilde{\Omega}$, an isomorphism $U\colon L^p(\widetilde{\Omega})\to 
L^p(\widetilde{\Omega})$ such that $\{U^n\, :\, n \in\Zdb\}$ is bounded, as well 
as two bounded maps $L^p(\Omega)\mathop{\longrightarrow}\limits^{J} L^p(\widetilde{\Omega})
\mathop{\longrightarrow}\limits^{Q}L^p(\Omega)$ such that $T^n=QU^nJ$ for any $n\geq 0$.
We also investigate functional calculus properties of Ritt operators and analogs of the above
results on noncommutative $L^p$-spaces.
\end{abstract}

\maketitle


\bigskip

\section{Introduction}
Let $(\Omega,\mu)$ be a measure space and let $1<p<\infty$. For any bounded
operator $T\colon L^p(\Omega)\to L^p(\Omega)$,  
consider the `square function'
\begin{equation}\label{1SF}
\norm{x}_{T,1}\,=\, \biggnorm{\biggl(\sum_{k=1}^{\infty} k\bigl 
\vert T^{k}(x) - T^{k-1}(x)\bigr\vert^2\biggr)^{\frac{1}{2}}}_{L^p},
\end{equation}
defined for any $x\in L^p(\Omega)$. 
Such quantities frequently appear in the analysis of $L^p$-operators.
They go back at least to  \cite{S2}, where they were used
in connection with martingale square functions to study 
diffusion semigroups and their discrete counterparts. 
Similar square functions for continuous semigroups
played a key role in the recent development of 
$H^\infty$-calculus and its applications. See in particular 
the fundamental paper \cite{CDMY}, the survey \cite{LM3} and the references
therein.

It is shown in \cite{LMX1} that if $T$ is both a positive contraction and
a Ritt operator, then it satisfies a uniform estimate
$\norm{x}_{T,1}\lesssim\norm{x}_{L^p}$ for $x\in L^p(\Omega)$. This estimate and
related ones lead to strong maximal inequalities for this
class of operators (see also \cite{LMX2}). Next in the paper
\cite{LM0}, the second named author studies the operators $T$ such that both
$T\colon L^p(\Omega)\to L^p(\Omega)$ 
and its adjoint operator
$T^*\colon L^{p'}(\Omega)\to L^{p'}(\Omega)$ satisfy uniform
estimates
\begin{equation}\label{1DoubleSFE}
\norm{x}_{T,1}\,\lesssim\,\norm{x}_{L^p}\qquad\hbox{and}\qquad
\norm{y}_{T^*,1}\,\lesssim\,\norm{y}_{L^{p'}}
\end{equation}
for $x\in L^p(\Omega)$ and $y\in L^{p'}(\Omega)$. (Here $p'=\frac{p}{p-1}$ is the conjugate number of $p$.)
It is shown that 
(\ref{1DoubleSFE}) implies that $T$ is an $R$-Ritt operator (see Section 2 below
for the definition) and that (\ref{1DoubleSFE}) is equivalent to 
$T$ having a bounded $H^\infty$-calculus with respect to a Stolz domain of the 
unit disc with vertex at $1$.

The present paper is a continuation of these investigations. Our main result
is a characterization of (\ref{1DoubleSFE}) in terms of dilations. 
We show that (\ref{1DoubleSFE}) holds true if and only if $T$ is $R$-Ritt
and there exist another measure space $(\widetilde{\Omega},\widetilde{\mu})$, two bounded maps
$J\colon L^p(\Omega)\to L^p(\widetilde{\Omega})$ and
$Q\colon L^p(\widetilde{\Omega})\to L^p(\Omega)$, as well as an isomorphism
$U\colon L^p(\widetilde{\Omega})\to L^p(\widetilde{\Omega})$ such that
$\bigl\{U^n\, :\, n\in\Zdb\bigr\}$ is bounded 
and 
$$
T^n=QU^nJ,\qquad n\geq 0.
$$
This result will be established in Section 4. 
It should be regarded as a discrete analog of the main result
of \cite{FW}. 

In Section 3, we consider variants of (\ref{1SF}) as follows. 
Assume that $T\colon L^p(\Omega)\to L^p(\Omega)$ is a Ritt operator. Then
$I-T$ is a sectorial operator and one can define its fractional power
$(I-T)^\alpha$ for any $\alpha>0$. Then we consider
\begin{equation}\label{1Alpha}
\norm{x}_{T,\alpha}\,=\, \biggnorm{\biggl(\sum_{k=1}^{\infty} k^{2\alpha -1}\bigl 
\vert T^{k-1}(I-T)^\alpha x\bigr\vert^2\biggr)^{\frac{1}{2}}}_{L^p}
\end{equation}
for any $x\in L^p(\Omega)$. Our second main result (Theorem \ref{3Equiv} below) is that
when $T$ is an $R$-Ritt operator, then the square functions $\norm{\ }_{T,\alpha}$
are pairwise equivalent. This result of independent 
interest should be regarded as a discrete analog of \cite[Thm. 1.1]{LM1}. We prove it here as
it is a key step in our characterization of (\ref{1DoubleSFE}) in terms of dilations. 

Section 2 mostly contains preliminary results. Section 5 is devoted to complements on
$L^p$-operators and their functional calculus properties, in connection with 
$p$-completely bounded maps. 
Finally Section 6 contains generalizations to operators $T\colon X\to X$ 
on general Banach spaces $X$. We pay a 
special attention to noncommutative $L^p$-spaces, 
in the spirit of \cite{JLX}.

\bigskip
We end this introduction with a few notation. If $X$ is a Banach space, we let 
$B(X)$ denote the algebra of all bounded operators on $X$ and we let 
$I_X$ denote the identity operator on $X$ (or simply $I$ if there is no ambiguity
on $X$). For any $T\in B(X)$, we let 
$\sigma(T)$ denote the spectrum of $T$. If $\lambda\in\Cdb\setminus\sigma(T)$ (the resolvent 
set of $T$), we let $R(\lambda,T)=(\lambda I_X-T)^{-1}$ denote the corresponding 
resolvent operator. 
We refer the reader to \cite{Die} for general information on Banach space geometry.
We will frequently use Bochner spaces $L^p(\Omega;X)$, for which we refer to \cite{DU}.

For any $a\in\Cdb$ and $r>0$, we let $D(a,r)=\{z\in\Cdb\, :\, \vert z-a\vert<r\}$
and we let 
$\Ddb=D(0,1)$ denote the open unit disc centered at $0$. Also we let 
$\Tdb=\{z\in\Cdb\,:\,\vert z\vert =1\}$ denote its boundary.

Whenever $\O\subset\Cdb$ is a non empty open set, we let
$H^{\infty}(\O)$ denote the space of all bounded 
holomorphic functions $f\colon\O\to\Cdb$. This is a Banach
algebra for the norm
$$
\norm{f}_{H^{\infty}(\footnotesize{\O})}\,=\,\sup\bigl\{\vert f(z)\vert\, :\, z\in\O\bigr\}.
$$
Also we let $\P$ denote the algebra of all complex polynomials. 

In the above presentation and later on in the paper we will use $\lesssim$ to indicate
an inequality up to a constant which does not depend on the particular element to 
which it applies. Then $A(x)\approx B(x)$ will mean that we both have
$A(x)\lesssim B(x)$ and $B(x)\lesssim A(x)$.

\section{Preliminaries on $R$-boundedness and Ritt operators}
This section is devoted to definitions and preliminary results
involving $R$-boundedness (and the companion notion of $\gamma$-boundedness),
matrix estimates and Ritt operators. We deal with operators acting
on an arbitrary Banach space $X$ (as opposed to the next two sections,
where $X$ will be an $L^p$-space).

Let $(\varepsilon_k)_{k\geq 1}$ be a sequence of independent 
Rademacher variables on some probability
space $\Omega_0$. We let $\Rad(X)\subset L^2(\Omega_0;X)$ be the
closure of ${\rm Span}\bigl\{\varepsilon_k\otimes x\, :\, k\geq 1,\
x\in X\bigr\}$ in the Bochner space $L^2(\Omega_{0};X)$. 
Thus for any finite family $x_1,\ldots,x_n$ in $X$, we have
$$
\biggnorm{\sum_{k=1}^{n} \varepsilon_k\otimes x_k}_{\Rad(X)} \,=\,
\biggr(\int_{\Omega_0}\biggnorm{\sum_{k=1}^{n} \varepsilon_k(\omega)\,
x_k}_{X}^{2}\,d\omega\,\biggr)^{\frac{1}{2}}.
$$
We say that a set $F\subset B(X)$ is $R$-bounded provided that
there is a constant $C\geq 0$ such that for any
finite families $T_1,\ldots, T_n$ in $F$ and  $x_1,\ldots,x_n$
in $X$, we have
$$
\biggnorm{\sum_{k=1}^{n} \varepsilon_k\otimes T_k
(x_k)}_{\Rad(X)}\,\leq\, C\, \biggnorm{\sum_{k=1}^{n}
\varepsilon_k\otimes x_k}_{\Rad(X)}.
$$
In this case we let $R(F)$ denote the smallest possible $C$,
which is called the $R$-bound of $F$.

Let $(g_k)_{k\geq 1}$ denote a sequence of independent
complex valued, standard Gaussian
random variables on some probability space $\Omega_1$, and let 
$\Gauss(X)\subset L^2(\Omega_1;X)$ be the
closure of ${\rm Span}\bigl\{g_k\otimes x\, :\, k\geq 1,\
x\in X\bigr\}$. Then replacing the $\varepsilon_k$'s and $\Rad(X)$ by
the $g_k$'s and $\Gauss(X)$ in the above paragraph, we obtain the similar 
notion of $\gamma$-bounded set. The corresponding $\gamma$-bound of a set $F$
is denoted by $\gamma(F)$. 

These two notions are very close to each other, 
however we need to work with both of them in this paper. Comparing them,
we recall that any $R$-bounded set $F\subset B(X)$ is automatically
$\gamma$-bounded, with $\gamma(F)\leq R(F)$. Moreover if $X$ has a finite cotype,
then the Rademacher averages and the Gaussian averages are equivalent on $X$
(see e.g. \cite[Prop. 12.11 and Thm. 12.27]{Die}), hence $F$ is $R$-bounded if 
(and only if) it is $\gamma$-bounded.

$R$-boundedness was introduced in \cite{BG} and then developed in the fundamental paper
\cite{CPSW}. We refer to the latter paper and to \cite[Section 2]{KW} for a detailed presentation.
We recall two facts which are highly relevant for our paper.
First, the closure of the absolute convex hull of any $R$-bounded set is
$R$-bounded \cite[Lem. 3.2]{CPSW}. This implies the following.

\begin{lemma}\label{2Aco} Let $F\subset B(X)$ be an $R$-bounded set, let
$J\subset\Rdb$ be an interval and let $C\geq 0$ be a constant. Then the set
$$
\biggl\{\int_J a(t)V(t)\, dt\ \Big\vert\, V\colon J\to F\ \hbox{is continuous},
\ a\in L^1(J),\ \norm{a}_{L^1(J)}\leq C\,
\biggr\}
$$
is $R$-bounded.
\end{lemma}

Second, if $X=L^p(\Omega)$ is an $L^p$-space with $1\leq p<\infty$,
then $X$ has a finite cotype and we have an equivalence
\begin{equation}\label{Rad(Lp)}
\biggnorm{\sum_{k} \varepsilon_k\otimes x_k}_{\Rad(L^p(\Omega))}\,
\approx \, \biggnorm{\biggl(\sum_{k}
\vert x_k\vert^2\biggr)^{\frac{1}{2}}}_{L^p(\Omega)}
\end{equation}
for finite families $(x_k)_k$ of $L^p(\Omega)$. 
Consequently a set $F\subset B\bigl(L^p(\Omega)\bigr)$ is $R$-bounded
if and only if it is $\gamma$-bounded, if and only if there exists a constant 
$C\geq 0$ such that for any
finite families $T_1,\ldots, T_n$ in $F$ and $x_1,\ldots,x_n$
in $L^p(\Omega)$, we have
$$
\biggnorm{\biggl(\sum_{k=1}^{n}\bigl\vert T_k(x_k)\bigr\vert^2\biggr)^{\frac{1}{2}}}_{L^p(\Omega)}\,
\leq\, C\, 
\biggnorm{\biggl(\sum_{k=1}^{n}\bigl\vert x_k\bigr\vert^2\biggr)^{\frac{1}{2}}}_{L^p(\Omega)}\,.
$$

\bigskip
In the sequel we represent any element of $B(\ell^2)$ by an infinite matrix $[c_{ij}]_{i,j\geq 1}$ in the usual
way. Likewise for any integer $n\geq 1$, we identify the algebra $M_n$ of all $n\times n$
matrices with the space of linear maps $\ell^2_n\to \ell^2_n$. Clearly an infinite matrix $[c_{ij}]_{i,j\geq 1}$
represents an element of $B(\ell^2)$ (in the sense that it is the matrix associated to a
bounded operator $\ell^2\to \ell^2$) if and only if
$$
\sup_{n\geq 1}\bignorm{[c_{ij}]_{1\leq i,j\leq n}}_{B(\ell^2_n)}\,<\infty\,.
$$
For any $[c_{ij}]_{1\leq i,j\leq n}$ in $M_n$, we set
$$
\bignorm{[c_{ij}]}_{\rm reg}\,=\,\bnorm{\bigr[\vert
c_{ij}\vert\,\bigr]}_{B(\ell^2_n)}.
$$
This is the so-called `regular norm' of the operator $[c_{ij}]\colon \ell^2_n\to \ell^2_n$.

\begin{lemma}\label{2Pisier}
For any matrix $[c_{ij}]$ in $M_n$, the following assertions are equivalent.
\begin{itemize}
\item [(i)] We have $\bnorm{[c_{ij}]}_{\rm reg}\leq 1$.
\item [(ii)] There exist two matrices
$[a_{ij}]$ and $[b_{ij}]$ in $M_n$ such that
$c_{ij}= a_{ij}b_{ij}$ for any $i,j=1,\ldots, n$, and we both have
$$
\sup_{1\leq i\leq n}\sum_{j=1}^{n}\vert a_{ij}\vert^2\,\leq 1\qquad\hbox{and}\qquad
\sup_{1\leq j\leq n}\sum_{i=1}^{n}\vert b_{ij}\vert^2\,\leq 1.
$$
\end{itemize}
\end{lemma}

The implication `(ii)$\Rightarrow$(i)' is an easy application of the Cauchy-Schwarz
inequality. The converse is due to Peller \cite[Section 3]{Pe} (see also \cite{Ak}).
We refer to \cite{P1} and \cite[Sect. 1.4]{P5} for more about this result 
and complements on regular norms.

The following result extends the boundedness of the Hilbert matrix (which corresponds to the case
$\beta=\gamma=\frac{1}{2}$). We thank \'Eric Ricard for his precious help in devising this proof.

\begin{proposition}\label{2Hilbert} Let $\beta,\gamma>0$ be two positive real numbers.
Then the infinite matrix
$$
\biggl[\frac{i^{\beta- \frac{1}{2}} j^{\gamma- \frac{1}{2}}}{(i+j)^{\beta+\gamma}}\biggr]_{i,j\geq 1}
$$
represents an element of $B(\ell^2)$.
\end{proposition}

\begin{proof} For any $i,j\geq 1$, set
$$
c_{ij} = \frac{i^{\beta- \frac{1}{2}} j^{\gamma- \frac{1}{2}}}{(i+j)^{\beta+\gamma}},\qquad
a_{ij} = c_{ij}^{\frac{1}{2}} \bigl(\tfrac{i}{j}\bigr)^{\frac{1}{4}},\qquad\hbox{and}\qquad
b_{ij} = c_{ij}^{\frac{1}{2}} \bigl(\tfrac{j}{i}\bigr)^{\frac{1}{4}}\,.
$$
Then $c_{ij}=a_{ij}b_{ij}$ for any $i,j\geq 1$, hence by the easy implication
of Lemma \ref{2Pisier}, it suffices to show that
\begin{equation}\label{2Hilbert1}
\sup_{i\geq 1}\sum_{j=1}^{\infty} \vert a_{ij}\vert^2\,<\infty\,\qquad\hbox{and}\qquad
\sup_{j\geq 1}\sum_{i=1}^{\infty}\vert b_{ij}\vert^2\,<\infty\,.
\end{equation}
Fix some $i\geq 1$. For any $j\geq 1$, we have
$$
\vert a_{ij}\vert^2 = c_{ij}\bigl(\tfrac{i}{j}\bigr)^{\frac{1}{2}}\,
=\,\frac{i^\beta j^{\gamma-1}}{(i+j)^{\beta+\gamma}}\,.
$$
Hence 
$$
\sum_{j=1}^{\infty} \vert a_{ij}\vert^2\,=\, i^\beta\biggl(\frac{1}{(i+1)^{\beta+\gamma}}\, + 
\sum_{j=2}^{\infty}\,\frac{1}{j^{1-\gamma}(i+j)^{\beta+\gamma}}\biggr)\,.
$$
Looking at the variations of the function $t\mapsto 1/(t^{1-\gamma}(i+t)^{\beta+\gamma})$ on
$(1,\infty)$, we immediately deduce that
$$
\sum_{j=1}^{\infty} \vert a_{ij}\vert^2\,\leq \, 1 + 2 \, i^\beta \int_{1}^{\infty} \frac{1}{t^{1-\gamma}(i+t)^{\beta+\gamma}}\, dt\,.
$$
Changing $t$ into $it$ in the latter integral, we deduce that
$$
\sum_{j=1}^{\infty} \vert a_{ij}\vert^2\,\leq \, 1 
+ 2 \int_{0}^{\infty}\frac{1}{t^{1-\gamma}(1+t)^{\beta+\gamma}}\, dt\,.
$$
This upper bound is finite and does not depend on $i$, which proves the first half of (\ref{2Hilbert1}). The 
proof of the second half is identical.
\end{proof}

We record the following elementary lemma for later use.

\begin{lemma}\label{2Positive} Let $[c_{ij}]_{i,j\geq 1}$ and
$[d_{ij}]_{i,j\geq 1}$ be infinite matrices of nonnegative real numbers,
such that $c_{ij} \leq d_{ij}$ for any $i,j\geq 1$. If the matrix
$[d_{ij}]_{i,j\geq 1}$ represents an element of $B(\ell^2)$, then
the same holds for  $[c_{ij}]_{i,j\geq 1}$.
\end{lemma}

We will need the following classical fact (see e.g. \cite[Cor. 12.17]{Die}).

\begin{lemma}\label{2Gauss} Let $X$ be a Banach space and let
$[b_{ij}]_{\substack{1\leq i\leq n\\ 1\leq j\leq m}}$ be an element of $M_{n,m}$. Then for
any $x_1,\ldots,x_m$ in $X$, we have
$$
\biggnorm{\sum_{i=1}^{n}\sum_{j=1}^{m} g_i \otimes b_{ij} x_j}_{{\rm Gauss}(X)} \leq\,
\bignorm{[b_{ij}]}_{B(\ell^2_m;\ell^2_n)} \biggnorm{\sum_{j=1}^m g_j \otimes x_j}_{{\rm
Gauss}(X)}.
$$
\end{lemma}

That result does not remain true if we replace Gaussian variables by
Rademacher variables and this defect is the main reason why it is
sometimes easier to deal with $\gamma$-boundedness than with
$R$-boundedness.

\begin{proposition}\label{2Reg}
Let $X$ be a Banach space, let $F= \big\{T_{ij}\, :\, i,j\geq
1\big\}$ be a $\gamma$-bounded family of operators on $X$, let $n\geq 1$
be an integer and let $[c_{ij}]_{1\leq i,j\leq n}$ be an element of
$M_n$. Then for any $x_1,\ldots,x_n$ in $X$, we have
$$
\biggnorm{\sum_{i,j=1}^{n} g_i\otimes c_{ij} T_{ij}(x_j)}_{{\rm
Gauss}(X)}\, \leq\, \gamma(F)\,\bignorm{[c_{ij}]}_{\rm
reg}\,\biggnorm{\sum_{j=1}^{n} g_j\otimes x_j}_{{\rm Gauss}(X)}.
$$
\end{proposition}

\begin{proof}
We can assume that $\bignorm{[c_{ij}]}_{\rm reg} \leq 1$. By Lemma
\ref{2Pisier}, we can write $c_{ij}=a_{ij}b_{ij}$ with
\begin{equation}
\label{2esti} \sup_{1 \leq i\leq n}\sum_{j=1}^{n}|a_{ij}|^2\leq 1 \ \
\ \ \ \text{and} \ \ \ \ \ \sup_{1 \leq j\leq
n}\sum_{i=1}^{n}|b_{ij}|^2\leq 1.
\end{equation}
Let $(g_{i,j})_{i,j\geq 1}$ be a doubly indexed family of
independent Gaussian variables. For any integers $1\leq i,j\leq  n$,
we define 
$$
A(i)=\left[\begin{array}{cccc} a_{i1} & a_{i2} & \hdots & a_{in}
\end{array}\right]
\qquad\hbox{and}\qquad 
B(j)=\left[\begin{array}{cccc} b_{1j} & b_{2j} & \hdots & b_{nj}
\end{array}\right]^{T}.
$$
Then we consider the two matrices
$$
A=\Diag\bigl(A(1),\ldots,A(n)\bigr)\,\in M_{n,n^2}
\qquad\hbox{and}\qquad
B=\Diag\bigl(B(1),\ldots,B(n)\bigr)\,\in M_{n^2,n}.
$$
Let $x_1,\ldots,x_n\in X$. Applying Lemma \ref{2Gauss} successively to $A$ and $B$, we then have
\begin{align*}
\bgnorm{\sum_{i,j=1}^{n} g_i \otimes c_{ij} T_{ij}(x_j)}_{{\rm
Gauss}(X)}
& \, = \, \bgnorm{\sum_{i,j=1}^{n} g_i \otimes a_{ij}b_{ij}T_{ij}(x_j)}_{{\rm Gauss}(X)}\\
& \,\leq \,\norm{A}\,\bgnorm{\sum_{i,j=1}^{n} g_{ij}\otimes b_{ij}T_{ij}(x_j)}_{{\rm Gauss}(X)}\\
& \,\leq \, \gamma(F)\,\norm{A}\,\bgnorm{\sum_{i,j=1}^{n} g_{ij}\otimes b_{ij} x_j}_{{\rm Gauss}(X)}\\
& \,\leq\,  \gamma(F)\,\norm{A}\norm{B}\,\bgnorm{\sum_{j=1}^{n} g_{j}\otimes  x_j}_{{\rm Gauss}(X)}.
\end{align*}
We have
$$
\norm{A}=\,\sup_{1\leq i\leq n}\bignorm{A(i)}_{M_{1,n}}\,=\,\sup_{1\leq i\leq n}
\biggl(\sum_{j=1}^{n}|a_{ij}|^2\biggr)^{\frac{1}{2}},
$$
hence $\norm{A}\leq 1$ by (\ref{2esti}). Likewise, we have $\norm{B}\leq 1$ hence the
above inequality yields the result.
\end{proof}

\bigskip
We now turn to Ritt operators, the key class of this paper, 
and recall some of their main features. Details and
complements can be found in \cite{Bl1, Bl2, LM0, Ly, NZ, N, V}.
We say that an operator $T\in B(X)$ is a Ritt operator if the two sets
\begin{equation}\label{2R1}
\{T^n\,:\, n\geq 0\}\qquad\hbox{and}\qquad \bigl\{n(T^n-T^{n-1})
\,:\, n\geq 1\bigr\}
\end{equation}
are bounded. This is equivalent to the spectral inclusion
\begin{equation}\label{2R2}
\sigma(T)\subset \overline{\Ddb}
\end{equation}
and the boundedness of the set
\begin{equation}\label{2R3}
\bigl\{(\lambda -1)R(\lambda,T)\, :\, \vert \lambda\vert>1\bigr\}.
\end{equation}
This resolvent estimate outside the unit disc is called the `Ritt condition'.

Likewise we say that $T$ is an $R$-Ritt operator if the two sets in (\ref{2R1}) are $R$-bounded. This
is equivalent to the inclusion (\ref{2R2}) and the $R$-boundedness of the set (\ref{2R3}).

For any angle $\gamma\in\bigl(0,\frac{\pi}{2}\bigr)$, let
$B_\gamma$ be the interior of the convex hull of $1$ and the 
disc $D(0,\sin\gamma)$ (see Figure 1 below).

\begin{figure}[ht]
\vspace*{2ex}
\begin{center}
\includegraphics[scale=0.4]{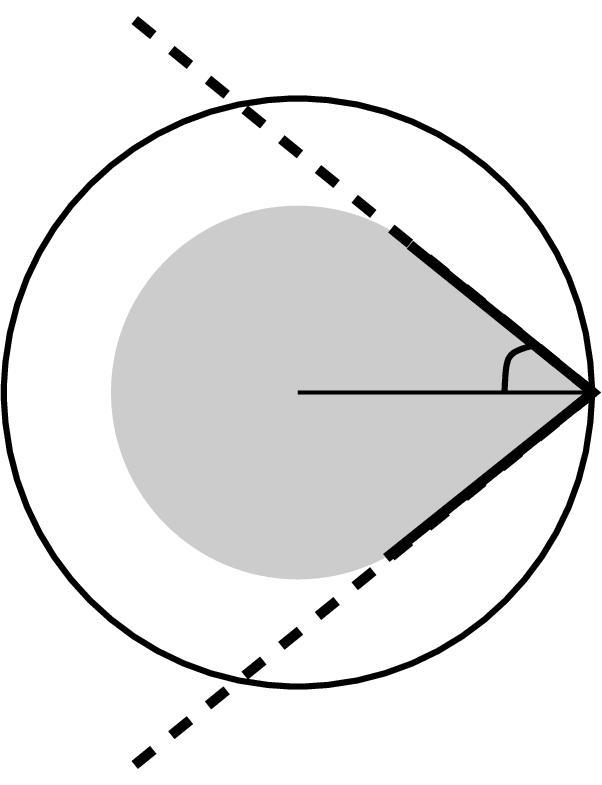}
\begin{picture}(0,0)
\put(-2,65){{\footnotesize $1$}}
\put(-68,65){{\footnotesize $0$}}
\put(-32,81){{\footnotesize $\gamma$}}
\put(-55,85){{\small $B_\gamma$}}
\end{picture}
\end{center}
\caption{\label{f1}}
\end{figure}

\noindent
Then the Ritt condition and its $R$-bounded version can be strengthened as follows.

\begin{lemma}\label{2Blunck} Let $T\colon X\to X$ be a Ritt operator
(resp. an $R$-Ritt operator). There exists an angle $\gamma\in (0,\frac{\pi}{2})$ such that
\begin{equation}\label{2Blunck1}
\sigma(T)\subset B_\gamma\cup\{1\}
\end{equation}
and the set
\begin{equation}\label{2Blunck2}
\bigl\{(\lambda-1)R(\lambda,T)\, :\, \lambda\in\Cdb\setminus 
B_\gamma,\ \lambda\not= 1\bigr\}
\end{equation}
is bounded (resp. $R$-bounded).
\end{lemma}

This essentially goes back to \cite{Bl1}, see \cite{LM0} for details.

For any angle $\theta\in(0,\pi)$, let 
\begin{equation}\label{2Sector}
\Sigma_\theta=\bigl\{z\in\Cdb\, :\, \bigl\vert{\rm Arg}(z)\bigr\vert <\theta\bigr\}
\end{equation}
be the open sector of angle $2\theta$ around the positive real axis $(0,\infty)$.
We say that a closed operator $A\colon D(A)\to X$ with dense domain $D(A)$ is sectorial
if there exists $\theta\in(0,\pi)$ such that $\sigma(A)\subset\overline{\Sigma_\theta}$
and the set 
\begin{equation}\label{2Sectorial}
\{zR(z,A)\, :\, z\in\Cdb\setminus \overline{\Sigma_\theta}\}
\end{equation}
is bounded.

Let $T$ be a Ritt operator and let $\gamma\in (0,\frac{\pi}{2})$ 
be such that the spectral inclusion
(\ref{2Blunck1}) holds true and the set (\ref{2Blunck2}) is bounded.
Then $A=I-T$ is a sectorial operator. Indeed $1-B_\gamma\subset \Sigma_\gamma$
and $zR(z,A)=\bigl((1-z)-1\bigr)R(1-z,T)$ for any $z\notin\overline{\Sigma_\gamma}$.
Hence for $\theta=\gamma$, the set (\ref{2Sectorial}) is bounded.
Thus for any $\alpha>0$, one can consider the fractional power $(I-T)^\alpha$. 
We refer e.g. to \cite[Chap. 3]{Ha} for various definitions of these 
(bounded) operators and their basic properties.
Fractional powers of Ritt operators can be expressed by a natural Dunford-Riesz functional
calculus formula. Indeed it was observed in \cite{LM0} that for any polynomial $\varphi$, we have
\begin{equation}\label{2DR}
\varphi(T)(I-T)^\alpha\,=\,\frac{1}{2\pi i }\,\int_{\partial B_\gamma} \varphi(\lambda)(1-\lambda)^\alpha
R(\lambda,T)\, d\lambda\,,
\end{equation}
where the countour $\partial B_\gamma$ is oriented counterclockwise.

P. Vitse proved in \cite{V} that if $T\colon X\to X$ is a Ritt operator, then
for any integer $N\geq 0$, the set $\{n^N T^{n-1}(I-T)^N\, :\, n\geq 1\}$ is bounded.
Our next statement is a continuation of these results.

\begin{proposition}\label{2Nev}
Let $X$ be a Banach space and let $T\colon X\to X$ be a Ritt operator
(resp. an $R$-Ritt operator).
For any $\alpha>0$, the set
$$
\bigl\{n^\alpha (rT)^{n-1}(I-rT)^{\alpha}\, :\, n\geq 1,\ r\in(0,1]\bigr\}
$$
is bounded (resp. $R$-bounded).
\end{proposition}

\begin{proof}
We will prove this result in the `$R$-Ritt case' only.
The `Ritt case' is similar and simpler. Assume that 
$T$ is $R$-Ritt. Applying Lemma \ref{2Blunck}, we let $\gamma\in (0,\frac{\pi}{2})$ be such that
(\ref{2Blunck1}) holds true and the set (\ref{2Blunck2}) is $R$-bounded.
Let $r\in(0,1]$ and let $\lambda\in\Cdb\setminus B_\gamma$, with $\lambda\not=1$. 
Then $\frac{\lambda}{r}\in\Cdb\setminus B_\gamma$ hence 
$\frac{\lambda}{r}$ belongs to the resolvent set of $T$ and we have
$$
(\lambda-1)R(\lambda,rT) =\,\frac{\lambda -1}{\lambda -r}\Bigl(\frac{\lambda}{r}-1\Bigr) R\Bigl(
\frac{\lambda}{r},T\Bigr).
$$
Since the set 
$$
\Bigl\{\frac{\lambda -1}{\lambda -r}\,:\, \lambda\in \Cdb\setminus 
B_\gamma,\ \lambda\not=1,\ r\in (0,1]\Bigr\}
$$
is bounded, it follows from the above formula that the set 
\begin{equation}\label{2Set}
\bigl\{(\lambda-1)R(\lambda,rT)\,:\, \lambda\in \Cdb\setminus B_\gamma,\
\lambda\not=1,\ r\in (0,1]\bigr\}
\end{equation}
is $R$-bounded.

The boundary $\partial B_\gamma$ is the juxtaposition
of the segment $\Gamma_+$ going from $1$ to $1-\cos(\gamma) e^{-i\gamma}$, of the
segment $\Gamma_-$ going from $1-\cos(\gamma) e^{i\gamma}$ to $1$
and of the curve $\Gamma_0$ going from $1-\cos(\gamma) e^{-i\gamma}$ to 
$1-\cos(\gamma) e^{i\gamma}$ counterclockwise along the circle of center $0$ and radius
$\sin\gamma$.

Consider a fixed number $\alpha>0$. For any integer $n\geq 1$ and any $r\in
(0,1]$, we have
$$
(rT)^{n-1}(I-rT)^{\alpha}\,=\,\frac{1}{2\pi i}\int_{\partial B_\gamma} \lambda^{n-1}
(1-\lambda)^{\alpha} R(\lambda,rT)\, d\lambda\,
$$
by applying (\ref{2DR}) to $rT$. Hence we may write
$$
n^\alpha(rT)^{n-1}(I-rT)^{\alpha}\,=\,\frac{-n^\alpha}{2\pi i}
\int_{\partial B_\gamma} \lambda^{n-1}(1-\lambda)^{\alpha-1}(\lambda-1)R(\lambda,rT)\, d\lambda\,.
$$
According to the $R$-boundedness of the set (\ref{2Set}) and Lemma \ref{2Aco}, it therefore
suffices to show that
the integrals
$$
I_n\,=\, n^\alpha\int_{\partial B_\gamma} \vert\lambda\vert^n \vert 1-\lambda\vert^{\alpha-1}
\,\vert d\lambda\vert
$$
are uniformly bounded (for $n$ varying in $\Ndb$). Let us decompose each of these
integrals as $I_n=I_{n,0} + I_{n,+} + I_{n,-}$, with 
$$ 
I_{n,0}=n^\alpha\int_{\Gamma_{0}}\cdots\,\vert d\lambda\vert\,,\quad
I_{n,+}=n^\alpha\int_{\Gamma_{+}}\cdots\,\vert d\lambda\vert\,,\quad\hbox{and}\quad
I_{n,-}=n^\alpha\int_{\Gamma_{-}}\cdots\,\vert d\lambda\vert\,.
$$

For $\lambda\in \Gamma_0$, we both have 
$$
\cos\gamma\leq \vert 1-\lambda\vert\leq  2\qquad\hbox{and}\qquad 
\vert\lambda\vert=\sin\gamma. 
$$
Since the sequence $\bigl(n^\alpha(\sin\gamma)^n\bigr)_{n\geq 1}$
is bounded, this readily implies that the sequence $(I_{n,0})_{n\geq 1}$ is bounded.

Let us now estimate $I_{n,+}$. For any $t\in[0,\cos\gamma]$,
we have $t^2\leq t\cos\gamma$ hence
$$
\bigl\vert 1 -te^{-i\gamma}\bigr\vert^2\, =\, 1+t^2-2t\cos\gamma\leq 1-t\cos\gamma.
$$
Hence
$$
I_{n,+} = n^\alpha\int_{0}^{\cos\gamma}
\bigl\vert 1 -te^{-i\gamma}\bigr\vert^n t^{\alpha-1}\, dt\ \leq
n^\alpha\int_{0}^{\cos\gamma}
\bigl(1 -t\cos\gamma\bigr)^{\frac{n}{2}} t^{\alpha-1}\, dt\,.
$$
Changing $t$ into $s=t\cos\gamma$ and using the inequality $1-s\leq e^{-s}$, we deduce that
$$
I_{n,+} \leq\,\frac{n^\alpha}{(\cos\gamma)^{\alpha}}\int_{0}^{\cos^2\gamma}
s^{\alpha-1}\,e^{-\frac{sn}{2}} ds\,.
$$
This yields (changing $s$ into $u=\frac{sn}{2}$) 
$$
I_{n,+} \leq\,\frac{2^\alpha}{(\cos\gamma)^{\alpha}}\int_{0}^{\infty} u^{\alpha-1}e^{-u}\, du\,.
$$
Thus the sequence $(I_{n,+})_{n\geq 1}$ is bounded. Since $I_{n,-}=I_{n,+}$, 
this completes the proof of the boundedness of $(I_n)_{n\geq 1}$.
\end{proof}

\medskip
\section{Equivalence of square functions}
Throughout the next two sections, we fix
a measure space $(\Omega,\mu)$ and a number $1<p<\infty$. We shall deal 
with operators acting on the Banach space $X=L^p(\Omega)$.
We start with a precise definition of (\ref{1SF}) and (\ref{1Alpha})
and a few comments.

Let $T\colon L^p(\Omega)\to L^p(\Omega)$ be a bounded operator and
let $x\in L^p(\Omega)$. Let us consider
$$
x_k=k^{\frac{1}{2}}\big(T^{k}(x) - T^{k-1}(x)\big)
$$
for any $k\geq 1$. If the sequence
$(x_k)_{k\geq 1}$ belongs to the space $L^p(\Omega;\ell^2)$, then $\norm{x}_{T,1}$ is defined
as the norm of $(x_k)_{k\geq 1}$ in that space. Otherwise,
we set $\norm{x}_{T,1}=\infty$. If $T$ is a Ritt operator, then
the quantities $\norm{x}_{T,\alpha}$ are defined in a similar manner
for any $\alpha>0$. In particular, $\norm{x}_{T,\alpha}$ can be infinite.

These square functions are natural discrete analogs of the square functions 
asociated to sectorial operators (see \cite{CDMY} and the survey paper \cite{LM3}).

Assume that $T$ is a Ritt operator. Then $T$ is power bounded hence by
the Mean Ergodic Theorem (see e.g. \cite[Subsection 2.1.1]{Kr}), we have a direct sum
decomposition
\begin{equation}\label{3ME}
L^p(\Omega)\,=\,{\rm Ker}(I-T)\oplus\overline{{\rm Ran}(I-T)},
\end{equation}
where ${\rm Ker}(\cdotp)$ and ${\rm Ran}(\cdotp)$ denote the kernel and the range, 
respectively.
For any $\alpha>0$, we have ${\rm Ker}\bigl((I-T)^\alpha\bigr)={\rm Ker}(I-T)$. This implies that 
\begin{equation}\label{3Null}
\norm{x}_{T,\alpha}=0\,\Longleftrightarrow\, x\in {\rm Ker}(I-T).
\end{equation}

Given any $\alpha>0$, a general question is to determine whether 
$\norm{x}_{T,\alpha}<\infty$ for any $x$ in $L^p(\Omega)$. It is easy to check,
using the Closed graph Theorem, that this finiteness property
is equivalent to the existence of a constant $C\geq 0$ such that
\begin{equation}\label{3SFE}
\norm{x}_{T,\alpha}\leq C\norm{x}_{L^p},\qquad x\in L^p(\Omega).
\end{equation}
In \cite{LM0}, the second named author established the 
following connection between the boundedness of discrete square functions and 
functional calculus properties.

\begin{theorem}\label{3CF} (\cite{LM0})
Let $T\colon L^p(\Omega)\to L^p(\Omega)$ be a Ritt operator, with $1<p<\infty$. 
The following assertions are equivalent.
\begin{enumerate}
\item [(i)] The operator $T$ and its adjoint $T^*\colon L^{p'}(\Omega)\to L^{p'}(\Omega)$
both  satisfy uniform estimates
$$
\norm{x}_{T,1}\,\lesssim\,\norm{x}_{L^p}\qquad\hbox{and}\qquad
\norm{y}_{T^*,1}\,\lesssim\,\norm{y}_{L^{p'}}
$$
for $x\in L^p(\Omega)$ and $y\in L^{p'}(\Omega)$.
\item [(ii)] There exists an angle 
$0<\gamma<\frac{\pi}{2}$ and a constant $K\geq 0$ such that
$$
\norm{\varphi(T)}\,\leq\, K  \,\norm{\varphi}_{H^{\infty}(B_\gamma)}
$$
for any $\varphi\in \P$.
\item [(iii)] The operator $T$ is $R$-Ritt and there exists an angle 
$0<\theta<\pi$ such that $I-T$ admits a bounded $H^\infty(\Sigma_\theta)$ functional calculus.
\end{enumerate}
\end{theorem}

Besides \cite{LM0},
we refer to \cite{CDMY, KW, LM4, MI}
for general information on $H^\infty(\Sigma_\theta)$ functional calculus for sectorial operators.

The main purpose of this section is to show that if $T$ is $R$-Ritt, then the square functions
$\norm{\ }_{T,\alpha}$ are pairwise equivalent. Thus the existence of an estimate (\ref{3SFE})
does not depend on $\alpha>0$. This result (Theorem \ref{3Equiv} below) is a discrete 
analog of the equivalence of square functions associated to $R$-sectorial operators,
as established in \cite{LM1}.

We start with preliminary results which allow to estimate square functions $\norm{x}_{T,\alpha}$
by means of approximation processes.

\begin{lemma}\label{3Approx}
Assume that $T\colon L^p(\Omega)\to L^p(\Omega)$ is a Ritt operator, and let $\alpha>0$.
\begin{enumerate}
\item [(1)] 
For any operator $V\colon L^p(\Omega)\to L^p(\Omega)$ such that $VT=TV$ and any $x\in L^p(\Omega)$, we have
$$
\norm{V(x)}_{T,\alpha}\leq\norm{V}\norm{x}_{T,\alpha}.
$$
\item [(2)] 
For any $x\in {\rm Ran}(I-T)$, we have
$\norm{x}_{T,\alpha}<\infty$.
\item [(3)] Let $\nu\geq\alpha+1$ be an integer and let $x\in {\rm Ran}\bigl((I-T)^\nu\bigr)$. Then
$$
\norm{x}_{T,\alpha}\,=\,\lim_{r\to 1^{-}}\norm{x}_{rT,\alpha}.
$$
\end{enumerate}
\end{lemma}
\begin{proof}
(1): Consider  $V\in B\bigl(L^p(\Omega)\bigr)$. As is well-known,
the tensor product $V\otimes I_{\ell^2}$ extends to a bounded operator
$V\overline{\otimes}I_{\ell^2}\colon L^p(\Omega;\ell^2)\longrightarrow L^p(\Omega;\ell^2)$, with
$\norm{V\overline{\otimes}I_{\ell^2}}=\norm{V}$. Assume that $VT=TV$ and let $x$ be such 
that $\norm{x}_{T,\alpha}<\infty$. Then 
we have
$$
\bigl(k^{\alpha-\frac{1}{2}}T^{k-1}(I-T)^\alpha\bigl(V(x)\bigr)\bigr)_{k\geq 1}\, =\, V\overline{\otimes}I_{\ell^2} \,\Bigl[\bigl(k^{\alpha-\frac{1}{2}}T^{k-1}(I-T)^\alpha(x)\bigr)_{k\geq 1}\Bigr],
$$
and the result follows at once.

\smallskip
(2): Assume that $x=(I-T)x'$ for some $x'\in L^p(\Omega)$. By Proposition
\ref{2Nev}, there exists a constant $C$ such that
\begin{align*}
 \sum_{k=1}^{\infty}\bnorm{k^{\alpha-\frac{1}{2}}T^{k-1}(I-T)^{\alpha}(x)}_{L^p}
&\,=\, \sum_{k=1}^{\infty}k^{\alpha-\frac{1}{2}}\bnorm{T^{k-1}(I-T)^{\alpha+1}(x')}_{L^p} \\
&\,\leq\, \norm{x'}_{L^p}\sum_{k=1}^{\infty}k^{\alpha-\frac{1}{2}}\,\frac{C}{k^{\alpha+1}} \\
&\,\leq \,  C\, \norm{x'}_{L^p}\sum_{k=1}^{\infty}{k^{-\frac{3}{2}}} \,<    \infty.
\end{align*}
This implies that $\bigl(k^{\alpha-\frac{1}{2}}T^{k-1}(I-T)^{\alpha}(x)\bigr)_{k\geq 1}$
belongs to $L^p(\Omega;\ell^2)$.

\smallskip
(3): It is clear that 
$(I-rT)^\alpha\to (I - T)^\alpha$ when $r\to 1^{-}$. 
Assume that $x\in {\rm Ran}\bigl((I-T)^\nu\bigr)$. Arguing as in part (2) we find that
the sequence $\bigl(k^{\alpha-\frac{1}{2}}T^{k-1}(x)\bigr)_{k\geq 1}$
belongs to $L^p(\Omega;\ell^2)$. Then arguing as in part (1), we obtain that 
$$
\Bignorm{\Bigl(
k^{\alpha-\frac{1}{2}}T^{k-1}\bigl((I-rT)^\alpha -(I-T)^\alpha\bigr)
(x)\Bigr)_{k\geq 1}}_{L^p(\ell^2)}\longrightarrow\, 0
$$
when $r\to 1^{-}$. This implies the convergence result.
\end{proof}

\begin{theorem}\label{3Equiv}
Assume that $T\colon L^p(\Omega)\to L^p(\Omega)$ is an $R$-Ritt
operator. Then for any $\alpha,\beta>0$, we have an equivalence
$$
\norm{x}_{T,\alpha}\,\approx\,\norm{x}_{T,\beta},\qquad x\in L^p(\Omega).
$$
\end{theorem}

\begin{proof}
We fix $\gamma>0$ such that $\alpha+\gamma$ is an integer $N\geq 1$.
For any integer $k\geq 1$, we define the complex number
$$
c_k=\frac{k(k+1)\cdots (k+N-2)}{k^{\alpha-\frac{1}{2}}},
$$
with the convention that $c_k=\frac{1}{k^{\alpha-\frac{1}{2}}}$ if $N=1$.
For any $z\in\Ddb$, we have
$$
\sum_{k=1}^{\infty}k(k+1)\cdots
(k+N-2)z^{k-1}\,=\,\frac{(N-1)!}{(1-z)^N}.
$$
Hence
$$
\sum_{k=1}^{\infty}c_kk^{\alpha-\frac{1}{2}}z^{2k-2}(1-z^2)^N
\,=\,\sum_{k=1}^{\infty}k(k+1)\cdots (k+N-2)(z^2)^{k-1}(1-z^2)^N
\,=\, (N-1)!.
$$
Since the operator $T$ is power bounded, we deduce that for every $r\in(0,1)$ we have
$$
\sum_{k=1}^{\infty}c_kk^{\alpha-\frac{1}{2}}(rT)^{2k-2}\big(I-(rT)^2\big)^N\,=\,(N-1)!I,
$$
the series in the left handside being normally convergent.
Since $(I+rT)^{N}$ is invertible, this yields
$$
\sum_{k=1}^{\infty}c_k(rT)^{k-1}(I-rT)^{\gamma}k^{\alpha-\frac{1}{2}}(rT)^{k-1}
(I-rT)^\alpha=(N-1)!(I+rT)^{-N}.
$$
Let $x\in L^p(\Omega)$. For any integer $m\geq 1$ and any $r\in(0,1)$, we let
$$
y_m(r)\,=\, (N-1)!(I+rT)^{-N}m^{\beta-\frac{1}{2}}(rT)^{m-1}(I-rT)^\beta x.
$$
Then it follows from the above identity that
$$
y_m(r)\,=\, 
\sum_{k=1}^{\infty}  c_k m^{\beta-\frac{1}{2}}(rT)^{m+k-2}(I-rT)^{\beta+\gamma}
\cdot k^{\alpha-\frac{1}{2}}(rT)^{k-1}(I-rT)^\alpha x.
$$
For any $n\geq 1$, we consider the partial sum
$$
y_{m,n}(r)\,=\, 
\sum_{k=1}^{n}  c_k m^{\beta-\frac{1}{2}}(rT)^{m+k-2}(I-rT)^{\beta+\gamma}
\cdot k^{\alpha-\frac{1}{2}}(rT)^{k-1}(I-rT)^\alpha x,
$$
and we have $y_{m,n}(r)\to y_{m}(r)$ when $n\to\infty$.

Let us write 
\begin{equation}\label{3Decomp}
c_k m^{\beta-\frac{1}{2}}  (rT)^{m+k-2}(I-rT)^{\beta+\gamma} \,=\,
\frac{m^{\beta-\frac{1}{2}} c_k }{(m+k-1)^{\gamma+\beta}}\bigl[(m+k-1)^{\gamma+\beta}(rT)^{m+k-2}(I-rT)^{\beta+\gamma}
\bigr]
\end{equation}
for any $m,k\geq 1$. Since $c_k \sim_{+\infty} k^{\gamma-\frac{1}{2}}$, there
exists a positive constant $K$ such that
$$
\frac{m^{\beta-\frac{1}{2}} c_k }{(m+k-1)^{\gamma+\beta}} \leq K\,
\frac{m^{\beta-\frac{1}{2}}k^{\gamma-\frac{1}{2}}}{(m+k)^{\gamma+\beta}}
$$
for any $m,k\geq 1$.
It therefore follows from Proposition \ref{2Hilbert} and Lemma \ref{2Positive} 
that the matrix
$$
\biggl[\frac{m^{\beta-\frac{1}{2}} c_k}{(m+k-1)^{\gamma+\beta}}\biggr]_{m,k\geq 1}
$$
represents an element of $B(\ell^2)$. Moreover, by Proposition
\ref{2Nev}, the set
$$
F\,=\,\bigl\{ (m+k-1)^{\gamma+\beta}(rT)^{m+k-2}(I-rT)^{\gamma+\beta}\,:\,
m,k\geq 1,\ r\in (0,1] \bigr\}
$$
is $R$-bounded. Hence by (\ref{Rad(Lp)}), (\ref{3Decomp}) and  Proposition \ref{2Reg}, we 
get an estimate
$$
\biggnorm{\biggl(\sum_{m=1}^{M} \bigl\vert
y_{m,n}(r)\bigr\vert^2 \biggr)^{\frac{1}{2}}}_{L^p}\,\lesssim\,
\biggnorm{\biggl(\sum_{k=1}^{\infty}k^{2\alpha-1}
\bigl\vert (rT)^{k-1}(I-rT)^{\alpha}x\bigr\vert^2\biggr)^{\frac{1}{2}}}_{L^p}
$$
for any integer $M\geq 1$. Passing to the limit, we deduce that
$$
\biggnorm{\biggl(\sum_{m=1}^{\infty} \bigl\vert
y_{m}(r)\bigr\vert^2 \biggr)^{\frac{1}{2}}}_{L^p}\,\lesssim\,
\biggnorm{\biggl(\sum_{k=1}^{\infty}k^{2\alpha-1}
\bigl\vert (rT)^{k-1}(I-rT)^{\alpha}x\bigr\vert^2\biggr)^{\frac{1}{2}}}_{L^p}\,.
$$
Since the set 
$\bigl\{(N-1)!^{-1}(I+rT)^{N}\, :\, r\in (0,1) \bigr\}$ is bounded, we finally obtain
that
$$
\norm{x}_{rT,\beta}\,\lesssim\,\norm{x}_{rT,\alpha}.
$$
It is crucial to note that in this estimate, the majorizing constant hidden in the symbol
$\lesssim$ does not depend on $r\in(0,1)$.

Now let $\nu$ be an integer such that $\nu\geq \alpha +1$ and 
$\nu\geq \beta +1$.
Applying Lemma \ref{3Approx} (3), we deduce a uniform estimate
$$
\norm{x}_{T,\beta} \lesssim \norm{x}_{T,\alpha}
$$
for $x\in {\rm Ran}\bigl((I - T)^\nu\bigr)$. Next for any integer $m\geq 0$, set
$$
\Lambda_m=\frac{1}{m+1}\sum_{k=0}^{m}(I-T^k).
$$
It is clear that $\Lambda_m^\nu$ maps $L^p(\Omega)$ into ${\rm Ran}\bigl((I - T)^\nu\bigr)$.
Hence we actually have a uniform estimate
$$
\bnorm{\Lambda_m^\nu(x)}_{T,\beta}\,\lesssim\,
\bnorm{\Lambda_m^\nu(x)}_{T,\alpha},\qquad x\in L^p(\Omega),\ m\geq 1.
$$
Since $T$ is power bounded, the sequence $(\Lambda_m)_{m\geq 0}$ is bounded.
Applying Lemma \ref{3Approx} (1), we deduce a further
uniform estimate
$$
\bnorm{\Lambda_m^\nu(x)}_{T,\beta}\,\lesssim\,
\norm{x}_{T,\alpha},\qquad x\in L^p(\Omega),\ m\geq 1.
$$
Equivalently, we have 
$$
\biggnorm{\bigg(\sum_{k=1}^{l}  k^{2\beta-1} \bigl\vert T^{k-1}(I-T)^\beta
\Lambda_m^\nu(x)\bigr\vert^2\bigg)^{\frac{1}{2}}}_{L^p}\,\lesssim\,
\norm{x}_{T,\alpha},\qquad x\in X,\ m\geq 1,\ l\geq 1.
$$
For any $x\in\overline{{\rm Ran}(I-T)}$, $\Lambda_m(x)\to x$ and hence 
$\Lambda_m^\nu(x)\to x$ when $m\to\infty$.
Hence passing to the limit in the above inequaliy, we obtain a uniform estimate
$\norm{x}_{T,\beta} \lesssim \norm{x}_{T,\alpha}$ 
for  $x$ in $\overline{{\rm Ran}(I-T)}$.
Switching the roles of $\alpha$ and $\beta$, this shows that
$\norm{\ }_{T,\beta}$ and $\norm{\ }_{T,\alpha}$ are equivalent on
the space $\overline{{\rm Ran}(I-T)}$. Moreover 
$\norm{\ }_{T,\beta}$ and $\norm{\ }_{T,\alpha}$ vanish on ${\rm Ker}(I-T)$
by (\ref{3Null}). Appealing to the direct sum decomposition (\ref{3ME}),
we finally obtain that $\norm{\ }_{T,\beta}$ and $\norm{\ }_{T,\alpha}$ are equivalent on
$L^p(\Omega)$.
\end{proof}

The techniques developed so far in this paper allow us to prove the following proposition,
which complements Theorem \ref{3CF}. For a Ritt operator $T$, we let $P_T$ 
denote the projection onto ${\rm Ker}(I-T)$ which vanishes on $\overline{{\rm Ran}(I-T)}$
(recall (\ref{3ME})).

\begin{proposition}\label{3DbSFE}
Let $T\colon L^p(\Omega)\to L^p(\Omega)$ be a Ritt operator, with $1<p<\infty$. Then condition (i)
in Theorem \ref{3CF} is equivalent to:
\begin{enumerate}
\item [(i)'] 
We have an equivalence
$$
\norm{x}_{L^p}\,\approx\,\norm{P_T(x)}_{L^p}\,+\,\norm{x}_{T,1}
$$
for $x\in L^p(\Omega)$.
\end{enumerate}
\end{proposition}

\begin{proof}
That (i) implies (i)' was proved in \cite[Rem. 3.4]{LMX1} in the case when $T$ is `contractively regular'.
The proof in our present case is the same.

Assume (i)'. Let $y\in L^{p'}(\Omega)$. We consider a finite sequence $(x_k)_{k\geq 1}$ in 
$L^p(\Omega)$ and we set
$$
x\,=\,\sum_{k} k^{\frac{1}{2}} T^{k-1}(I-T)x_k\,.
$$
Then
$$
\biggl\vert\sum_{k} 
\bigl\langle k^{\frac{1}{2}} (T^*)^{k-1}(I-T^*)y , x_k \bigr\rangle\biggr\vert\,
=\,\vert\langle y,x\rangle\vert\,\leq\,\norm{x}_{L^p}\norm{y}_{L^{p'}}.
$$
Moreover $x\in{\rm Ran}(I-T)$ hence applying (i)', we deduce
$$
\biggl\vert\sum_{k} 
\bigl\langle k^{\frac{1}{2}} (T^*)^{k-1}(I-T^*)y,x_k \bigr\rangle\biggr\vert\,
\lesssim\,\norm{y}_{L^{p'}}\norm{x}_{T,1}.
$$
We will now show an estimate
\begin{equation}\label{3DbSFE1}
\norm{x}_{T,1}\,\lesssim\,\biggnorm{\biggl(\sum_{k}\vert x_k\vert^2\biggr)^{\frac{1}{2}}}_{L^p}.
\end{equation}
Then passing to the supremum over all finite sequences $(x_k)_{k\geq 1}$ in the unit ball of $L^p(\Omega;\ell^2)$,
we deduce that $\norm{y}_{T^*,1}\lesssim\norm{y}_{L^{p'}}$.

To show (\ref{3DbSFE1}), first note that for any integer $m\geq 1$, we may write
$$
m^{\frac{1}{2}} T^{m-1}(I-T)x\,=\,\sum_k \frac{m^{\frac{1}{2}}k^{\frac{1}{2}}}{(m+k)^2}\, 
(m+k)^2 T^{m+k-2}(I-T)^2 x_k\,.
$$
Second according 
to \cite{LM0}, the assumption (i)' implies that $T$ is an $R$-Ritt operator.
By Proposition \ref{2Nev}, the set $\{(m+k-1)^2 T^{m+k-2}(I-T)^2\, :\, m,k\geq 1\}$, and hence 
the set
$$
\bigl\{(m+k)^2 T^{m+k-2}(I-T)^2\, :\, m,k\geq 1\bigr\}
$$
is $R$-bounded. Therefore applying Propositions \ref{2Hilbert} and \ref{2Reg} we obtain
(\ref{3DbSFE1}).
\end{proof}

\medskip
\section{Loose dilations}
We will focus on the following notion of dilation for $L^p$-operators.

\begin{definition}\label{4loose} Let $T\colon L^p(\Omega)\to L^p(\Omega)$
be a bounded operator.
We say that it admits a {\it loose dilation} if
there exist a measure space $(\widetilde{\Omega},\widetilde{\mu})$, two bounded maps
$J\colon L^p(\Omega)\to L^p(\widetilde{\Omega})$ and
$Q\colon L^p(\widetilde{\Omega})\to L^p(\Omega)$, as well as an isomorphism
$U\colon L^p(\widetilde{\Omega})\to L^p(\widetilde{\Omega})$ such that
$\bigl\{U^n\, :\, n\in\Zdb\bigr\}$ is bounded  and
$$
T^n=QU^nJ,\qquad n\geq 0.
$$
\end{definition}

That notion is strictly weaker than the following more 
classical one.

\begin{remark}\label{4Strict}
We say that a bounded operator $T\colon L^p(\Omega)\to L^p(\Omega)$
admits a {\it strict dilation} if there exist a measure space 
$(\widetilde{\Omega},\widetilde{\mu})$, two contractions
$J\colon L^p(\Omega)\to L^p(\widetilde{\Omega})$ and
$Q\colon L^p(\widetilde{\Omega})\to L^p(\Omega)$, as well as an isometric isomorphism
$U\colon L^p(\widetilde{\Omega})\to L^p(\widetilde{\Omega})$ such that
$T^n=QU^nJ$ for any $n\geq 0$. 

This strict dilation property implies that $T$ is a contraction and that 
$J$ and $Q^*$ are both isometries. 

Conversely in the case $p=2$, Nagy's dilation Theorem (see e.g. \cite[Chapter 1]{NF}) ensures that any contraction 
$L^2(\Omega)\to L^2(\Omega)$ admits a strict dilation. 

Next, assume that $1<p\not=2<\infty$. Then it follows from \cite{Ak, A, CRW, Pe} 
that $T\colon L^p(\Omega)\to L^p(\Omega)$
admits a strict dilation if and only if there exists a positive contraction
$S\colon L^p(\Omega)\to L^p(\Omega)$ such that $\vert T(x)\vert\leq S(\vert x\vert)$
for any $x\in L^p(\Omega)$.
\end{remark}

Except for $p=2$ (see Remark \ref{4p=2} below), there is no
similar description of operators admitting a loose dilation.
The general issue behind our investigation is to try to characterize
the $L^p$-operators which satisfy this property. 
Theorem \ref{4Dilatable} below gives a satisfactory answer for the
class of Ritt operators.

\begin{remark}\label{4p=2}
Let $H$ be a Hilbert space, let $T\colon H\to H$ be a bounded
operator and let us say that $T$ admits a loose dilation
if
there exist a Hilbert space $K$, two bounded maps
$J\colon H\to K$ and
$Q\colon K\to H$, and an isomorphism $U\colon K\to K$ such that
$\bigl\{U^n\, :\, n\in\Zdb\bigr\}$ is bounded  and
$T^n=QU^nJ$ for any $n\geq 0$. Then this property 
is equivalent to $T$ being similar to a contraction.

Indeed assume that there exists an isomorphism $V\in B(H)$
such that $V^{-1}TV$ is a contraction. By Nagy's dilation
Theorem, that contraction admits a unitary dilation. In other words,
there is a unitary $U$ on a Hilbert space $K$ containing $H$, such that
$(V^{-1}TV)^n=q U^n j$ for any $n\geq 0$, where $j\colon H\to K$
is the canonical inclusion and $q=j^*$ is the corresponding orthogonal
projection. We obtain the loose dilation
property of $T$ by taking $J=j V^{-1}$ and $Q=Vq$.

The converse uses the notion of complete polynomial boundedness, for
which we refer to \cite{Pa1,Pa2}. Assume that $T$ admits a loose
dilation. Using \cite[Cor. 9.4]{Pa2} and elementary arguments, 
we obtain that $T$ is completely polynomially bounded. Hence 
it is similar to a contraction by \cite[Cor. 3.5]{Pa1}.
\end{remark}

According to the above result, the rest of this section is significant
only in the case $1<p\not=2<\infty$.

Let $S\colon \ell^p_{\Zdb} \to \ell^p_{\Zdb}$ denote the natural shift
operator given by $S \bigl((t_k)_{k\in \Zdb}\bigr)
=\bigl((t_{k-1})_{k\in\Zdb}\bigr)$. 
For any $\varphi$ in $\P$ (the algebra of complex polynomials), we set
\begin{equation}\label{4Pnorm}
\norm{\varphi}_{p}=\,\bnorm{\varphi(S)}_{B(\ell^p_{\Zdb})}.
\end{equation}
We recall that if $\varphi$ is given by $\varphi(z)=\sum_{k\geq 0} d_k z^k\,$, then 
$\varphi(S)$ is the convolution operator (with respect to the group $\Zdb$)
associated to the sequence $(d_k)_{k\in\Zdb}$. Alternatively, $\varphi(S)$
is the Fourier multiplier associated to the restriction 
of $\varphi$ to the unitary group $\Tdb$. 
We refer the reader to \cite{EG} for some elementary 
background on Fourier multiplier theory.

Let us decompose $(0,\pi)$ dyadically into
the following family $(I_j)_{j\in \mathbb{Z}}$ of intervals:
$$
I_j =\left\{
\begin{array}{cl}
\Bigl[\pi-\frac{\pi}{2^{j+1}},\pi-\frac{\pi}{2^{j+2}}\Bigr)    &  \text{if}\      j\geq 0      \\
\bigl[2^{j-1}\pi,2^j\pi\bigr)                     &  \text{if}\      j<0.          \\
\end{array}\right.
$$
Then we denote by $\Delta_j$ the corresponding arcs of $\Tdb$:
$$
\Delta_j=\big\{e^{it} \, :\,  t \in -I_j\cup I_j  \ \big\}.
$$
We will use the following version of the Marcinkiewicz
multiplier theorem (see \cite[Thm. 4.3]{Bl1} and also \cite{EG}).

\begin{theorem}\label{4Marci} 
Let $1<p<\infty$. Let $\phi\in L^\infty(\Tdb)$ and assume that
$\phi$ has uniformly bounded variations over the
$(\Delta_j)_{j\in \mathbb{Z}}$. Then $\phi$ induces a bounded
Fourier multiplier $M_\phi\colon \ell^p_\Zdb\to \ell^p_\Zdb$ and we have
$$
\bnorm{M_\phi}_{B(\ell^p_\Zdb)}\, \leq C_{p}\,\Bigl(\norm{\phi}_{L^{\infty}(\Tdb)}
+\,
\sup\bigl\{\var(\phi,\Delta_j)\, :\, j\in\Zdb\bigr\}\Bigr),
$$
where $\var(\phi,\Delta_j)$ is the usual variation of $\phi$ over
$\Delta_j$ and the constant $C_p$ only depends on $p$.
\end{theorem}

For convenience, 
  \ref{5ppb} and Proposition \ref{4Fourier} below
are given for an arbitrary Banach space $X$, although
we are mostly interested in the case when $X$ is an $L^p$-space.

\begin{definition}\label{5ppb}
We say that a bounded operator $T\colon X\to X$
is $p$-polynomially bounded if there exists a constant $C\geq 1$ such that
\begin{equation}\label{4Ppb}
\bnorm{\varphi(T)}\leq C\norm{\varphi}_{p}
\end{equation}
for any complex polynomial $\varphi$.
\end{definition}

The following connection with dilations is well-known.

\begin{proposition}\label{4Transference} 
If $T\colon L^p(\Omega)\to L^p(\Omega)$ admits a loose dilation, then it is
$p$-polynomially bounded.
\end{proposition}

\begin{proof} Assume that $T$ satisfies the dilation property given by Definition
\ref{4loose}. Then for any $\varphi\in\P$, we have
$\varphi(T)=Q\varphi(U)J$, hence 
$$
\bnorm{\varphi(T)}\leq
\norm{Q}\norm{J}\bnorm{\varphi(U)}.
$$ 
Moreover by the transference principle (see
\cite[Thm. 2.4]{CW}), $\bnorm{\varphi(U)}\leq K^2\norm{\varphi}_p$,
where $K\geq 1$ is any constant such that $\norm{U^n}\leq K$ for any
integer $n$. This yields the result.
\end{proof}

We will see in Section 5 that the converse of this proposition does not hold true.

The above proof shows that if $T\colon L^p(\Omega)\to L^p(\Omega)$ admits a strict dilation, then
$\norm{\varphi(T)}\leq\norm{\varphi}_p$ for any $\varphi\in\P$, a very
classical fact. The famous Matsaev Conjecture asks whether this inequality
holds for any $L^p$-contraction $T$ (even those with no strict dilation). 
This was disproved very recently 
by Drury in the case $p=4$ \cite{Dr}. It is unclear whether 
there exists an $L^p$-contraction $T$ satisfying $\norm{\varphi(T)}
\leq\norm{\varphi}_p$ for any $\varphi\in\P$,
without admitting a strict dilation.

\begin{proposition}\label{4Fourier} Let
$T\colon X\to X$ be a $p$-polynomially bounded
operator. Then $I-T$ is sectorial and 
for any $\theta\in\big(\tfrac{\pi}{2},\pi\big)$, it
admits a bounded $H^{\infty}(\Sigma_\theta)$ functional
calculus.
\end{proposition}

\begin{proof} 
Since $T$ is $p$-polynomially bounded, it is power bounded
hence $\sigma(T)\subset\overline{\Ddb}$. We can thus define 
$\varphi(T)$ for any rational function with poles outside 
$\overline{\Ddb}$. Furthermore (\ref{4Ppb}) holds as well 
for such functions, by approximation.

We fix two numbers $\tfrac{\pi}{2}<\theta<\theta'<\pi$ 
and we let (see Figure 2):
$$
\Ddb_\theta = \,D\bigl(-i\cot(\theta),\tfrac{1}{\sin(\theta)}\bigr)\cup
D\bigl(i\cot(\theta),\tfrac{1}{\sin(\theta)}\bigr).
$$

\begin{figure}[ht]
\begin{center}
\includegraphics[scale=0.5]{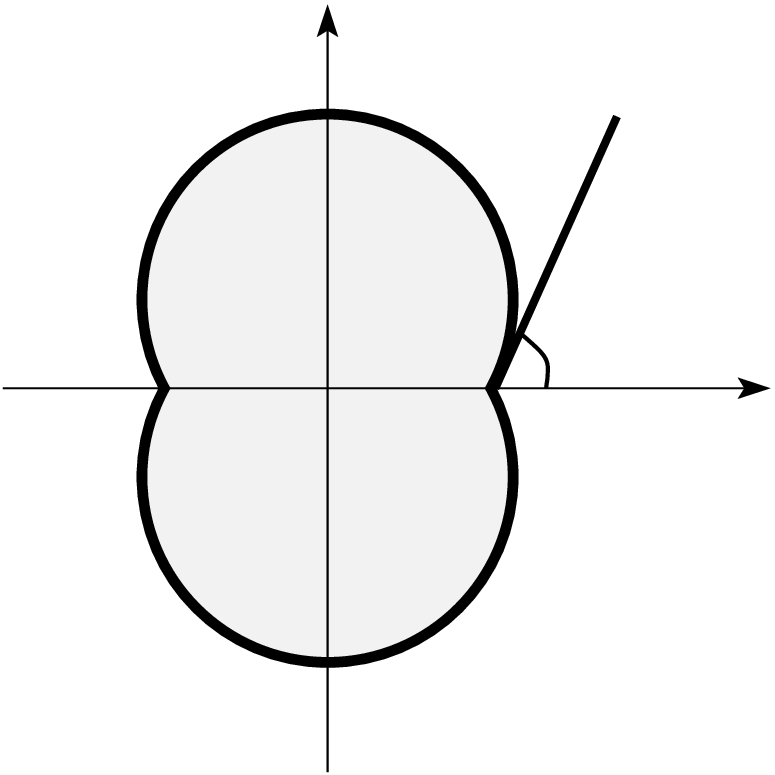}
\begin{small}
\begin{picture}(0,0) 
\put(-59,100){$\pi-\theta$} 
\put(-140,125){$\mathbb D_\theta$} 
\put(-155,83){$-1$} 
\put(-123,83){$0$} 
\put(-82,83){$1$} 
\end{picture}
\end{small}
\end{center}
\caption{\label{f2}}
\end{figure}

Clearly $\Ddb_\theta$ contains $\Ddb$.
For any $ t \in (-\pi,0)\cup(0,\pi)$, let $r(t)$ denote the radius
of the largest open disc centered at $e^{it}$ and included in $\Ddb_{\theta}$. 
If $t$ is positive and small enough, we have
\begin{align*}
 r(t)\,
 &= \,\frac{1}{\sin(\theta)}\, -\bigl\vert e^{it} +i\cot(\theta)\bigr\vert \\
 &= \,\frac{1}{\sin(\theta)}-\sqrt{\cos^2(t)+\big(\sin(t)+\cot(\theta)\big)^2}  \\
 &= \,\frac{1}{\sin(\theta)}\Bigl(1- \sqrt{1+2\sin(t)\sin(\theta)\cos(\theta)}\Bigr)\\
 &= -\cos(\theta)t+\frac{1}{2}\bigl(\sin(\theta)\cos^2(\theta)\bigr)t^2+O (t^3).
\end{align*}
Consequently, we have $r(t)>-\cos(\theta)t$ for $t>0$ small enough. We deduce that
if $j<0$ with $|j|$ large enough and $t \in I_{j}$, we have 
\begin{equation}\label{4Inclusion}
D\Bigl(e^{it},-\cos(\theta)\frac{\pi}{2^{|j|+1}}\Bigr)\subset 
\Ddb_{\theta}.
\end{equation}
The same holds for 
$t\in -I_j$. Moreover the intervals $I_j$ and $I_{-j}$
of the dyadic decomposition have length
equal to $\frac{\pi}{2^{|j|+1}}$.
Hence for any rational function $\varphi$ with poles outside $\overline{\Ddb_\theta}$
and any $j<0$ with $|j|$ large enough, we obtain that
\begin{align*}
   \var\bigl(\varphi_{\vert\Tdb} ,\Delta_j\bigr)\,
   & =\,  \int_{-I_j\cup I_j} \bigl\vert \varphi'(e^{it})\bigr\vert\, dt \\
   & \leq\, \int_{-I_j\cup I_j} \frac{\norm{\varphi}_{H^\infty(\Ddb_{\theta})}}{-\cos(\theta)\frac{\pi}{2^{|j|+1}}}\, dt \qquad \text{by (\ref{4Inclusion}) and Cauchy's inequalities,}\\
   & \leq\, \frac{\pi}{2^{|j|}}\,\cdotp\, \frac{\norm{\varphi}_{H^\infty(\Ddb_{\theta})}}{-\cos(\theta)\frac{\pi}{2^{|j|+1}}}
   \,=\,\frac{2\pi \norm{\varphi}_{H^\infty(\Ddb_{\theta})}}{-\cos(\theta)}\,.
\end{align*}
We have a similar result if $j \geq 0$ and large enough.
Applying Theorem \ref{4Marci}, we deduce a uniform estimate
$$
\bnorm{\varphi(S)}_{B(\ell^p_{\Zdb})}\, \lesssim\,
\norm{\varphi}_{H^\infty(\Ddb_{\theta})}.
$$
Combining with (\ref{4Ppb}) -as explained at the beginning of this proof- we obtain the existence of
a constant $K\geq 0$ such that for any rational function $\varphi$ with poles outside $\overline{\Ddb_\theta}$,
$$
\bnorm{\varphi(T)}_{B(X)}\, \leq\, K
\norm{\varphi}_{H^\infty(\Ddb_{\theta})}.
$$
Note that we have the following inclusion:
$$
1-\Ddb_\theta\subset\Sigma_\theta.
$$
Then let $\R_{\theta'}$ be the algebra of all rational functions 
with poles outside $\overline{\Sigma_{\theta'}}$ and with a nonpositive degree.
We deduce from above that for any $f\in \R_{\theta'}$,
$$
\bignorm{f(I- T)}  
\leq K\bignorm{f(1-\cdotp)}_{H^\infty(\Ddb_{\theta})} 
\leq K\norm{f}_{H^\infty(\Sigma_{\theta})}.
$$
According to \cite[Prop. 2.10]{LM4}, this readily implies that 
$I-T$ is sectorial and admits a bounded $H^\infty(\Sigma_{\theta'})$ 
functional calculus.
\end{proof}

\begin{theorem}\label{4Dilatable}
Let $T\colon L^p(\Omega)\to L^p(\Omega)$ be a 
Ritt operator, with $1<p<\infty$.
The following assertions are equivalent.
\begin{enumerate}
\item [(i)] The operator $T$ and its adjoint $T^*\colon L^{p'}(\Omega)\to L^{p'}(\Omega)$
both  satisfy uniform estimates
$$
\norm{x}_{T,1}\,\lesssim\,\norm{x}_{L^p}\qquad\hbox{and}\qquad
\norm{y}_{T^*,1}\,\lesssim\,\norm{y}_{L^{p'}}
$$
for $x\in L^p(\Omega)$ and $y\in L^{p'}(\Omega)$.
\item [(ii)] The operator $T$ is $R$-Ritt and admits a loose dilation.
\item [(iii)] The operator $T$ is $R$-Ritt and $p$-polynomially bounded.
\end{enumerate}
\end{theorem}

\begin{proof}
That (ii) implies (iii) follows from
Proposition \ref{4Transference}.

Assume (iii). By Proposition \ref{4Fourier}, $I-T$
admits a bounded $H^{\infty}(\Sigma_\theta)$ functional calculus for any $\theta>\tfrac{\pi}{2}$.
Since $T$ is R-Ritt, this implies (i) by Theorem \ref{3CF}.

Assume (i). It follows from \cite{LM0} (see Theorem \ref{3CF})
that $T$ is an $R$-Ritt operator.
Thus we only need to establish the dilation property of $T$.
Since $T$ is $R$-Ritt, Theorem \ref{3Equiv} ensures that
the square functions $\norm{\ }_{T,1}$ and $\norm{\ }_{T,\frac{1}{2}}$ are equivalent
on $L^{p}(\Omega)$.
Likewise, $\norm{\ }_{T^*,1}$ and $\norm{\ }_{T^*,\frac{1}{2}}$ are equivalent on
$L^{p'}(\Omega)$. Consequently assumption (i)
implies the existence of a constant $C\geq 1$ such that
$$
\norm{x}_{T,\frac{1}{2}}\,\leq\, C\norm{x}_{L^p}
\qquad\hbox{and}\qquad
\norm{y}_{T^*,\frac{1}{2}}\,\leq\, C\norm{y}_{L^{p'}}
$$
for any $x\in L^p(\Omega)$ and any $y\in L^{p'}(\Omega)$.

We will use the direct sum decomposition (\ref{3ME}), as well as the
analogous decompostion of $L^{p'}(\Omega)$ corresponding to $T^*$.
According to the above estimates, we may define two bounded maps
$$
j_1\colon \overline{{\rm Ran}(I-T)}\longrightarrow L^p(\Omega;\ell^2_\Zdb)
\qquad\hbox{and}\qquad
j_2\colon \overline{{\rm Ran}(I-T^*)}\longrightarrow L^{p'}(\Omega;\ell^2_\Zdb)
$$
as follows. For any $x\in \overline{{\rm Ran}(I-T)}$ and 
any $y\in \overline{{\rm Ran}(I-T^*)}$, 
we set 
$$
x_k=T^{k-1}(I-T)^{\frac{1}{2}}x
\qquad\hbox{and}\qquad
y_k=(T^*)^{k-1}(I-T^*)^{\frac{1}{2}}y
$$
if $k\geq 0$, we set 
$x_k=0$ and $y_k=0$ if $k<0$. Then we set
$$
j_1(x)=\,(x_k)_{k\in\Zdb}
\qquad\hbox{and}\qquad
j_2(y)=\,(y_k)_{k\in\Zdb}.
$$
Then we let $J_1\colon L^p(\Omega)\to L^p(\Omega)\mathop{\oplus}\limits^{p} L^p(\Omega;\ell^2_\Zdb)$
be the linear map taking any $x\in {\rm Ker}(I-T)$ to $(x,0)$ and any 
$x\in \overline{{\rm Ran}(I-T)}$ to $(0,j_1(x))$. We define 
$J_2\colon L^{p'}(\Omega) \to L^{p'}(\Omega)\mathop{\oplus}\limits^{p'} L^{p'}(\Omega;\ell^2_\Zdb)$
in a similar way.

For any $x\in \overline{{\rm Ran}(I-T)}$ and 
$y\in \overline{{\rm Ran}(I-T^*)}$, we have
\begin{align*}
\langle J_1 x, J_2 y  \rangle\, 
& =\,\sum_{k=1}^{\infty}\bigl\langle 
T^{k-1}(I-T)^{\frac{1}{2}}x,  (T^*)^{k-1}(I-T^*)^{\frac{1}{2}}y \bigr\rangle\\
& =\,\sum_{k=1}^{\infty}\bigl\langle 
T^{2(k-1)}(I-T)x, y \bigr\rangle\\
& =\,\sum_{k=1}^{\infty}\bigl\langle 
T^{2(k-1)}(I-T^2) (I+T)^{-1} x, y \bigr\rangle.
\end{align*}
For any integer $N\geq 1$
$$
\sum_{k=1}^{N} T^{2(k-1)}(I-T^2) \,=\, I -T^{2N}.
$$
Furthermore, $(I+T)^{-1} x$ belongs to $\overline{{\rm Ran}(I-T)}$
and the sequence $(T^{n})_{n\geq 0}$ strongly converges to $0$ on 
that subspace of $L^p(\Omega)$. Hence 
$$
\langle J_1 x, J_2 y\rangle \,=\, \langle 
(I+T)^{-1} x, y \rangle.
$$
Let $\Theta\colon L^p(\Omega)\to L^p(\Omega)$ be the linear map
taking any  $x\in \overline{{\rm Ran}(I-T)}$ to $(I+T)x$ and any 
$x\in {\rm Ker}(I-T)$ to itself. Then it follows from the above calculation that 
\begin{equation}\label{4Dil2}
\Theta J_2^* J_1\,=\,I_{L^p(\Omega)}.
\end{equation}
Let 
$$
Z=L^p(\Omega)\mathop{\oplus}\limits^{p} L^p(\Omega;\ell^2_\Zdb),
$$
and let $U\colon Z\to Z$ be the linear map which takes any $x\in L^p(\Omega)$ to itself and
any sequence $(x_k)_{k\in\Zdb}$ in $L^p(\Omega;\ell^2_\Zdb)$ to 
the shifted sequence $(x_{k+1})_{k\in\Zdb}$. Next let $P\colon Z\to Z$ be the
linear map which takes any $x\in L^p(\Omega)$ to itself and
any sequence $(x_k)_{k\in\Zdb}$ in $L^p(\Omega;\ell^2_\Zdb)$ to the truncated
sequence $(\ldots,0,\ldots,0,x_0,x_1,\ldots, x_k,\ldots)$. By construction, we have
\begin{equation}\label{4Dil3}
PU^nJ_1=J_1T^n,\qquad n\geq 0.
\end{equation}
We also have $J_2^*P=J_2^*$ hence setting $J=J_1\colon L^p(\Omega)\to Z$ and
$Q=\Theta J_2^*\colon Z\to L^p(\Omega)$, we deduce from (\ref{4Dil2}) and
(\ref{4Dil3}) that $T^n=QU^nJ$ for any $n\geq 0$. Furthermore, $U$ is an
isometric isomorphism on $Z$.
Thus we have established that $T$ satisfies the dilation property stated
in Definition \ref{4loose},
except that the dilation space is $Z$ instead of being an $L^p$-space. 

It is easy to modify the construction to obtain a dilation
through an $L^p$-space, as follows. First recall that 
using for example Gaussian variables, one can
isometrically represent $\ell^2_\Zdb$ as a complemented
subspace of an $L^p$-space (see e.g. \cite[Chapter 5]{P0}). The space 
$Z$ can be therefore represented as well as a complemented
subspace of an $L^p$-space. Thus we have 
$$
Z\oplus W\,=\, L^p(\widetilde{\Omega})
$$
for an appropriate measure space $(\widetilde{\Omega},\widetilde{\mu})$ and
some Banach space $W$.
Let $J'\colon L^p(\Omega)\to L^p(\widetilde{\Omega})$ be defined by 
$J'(x)=(J(x),0)$, let $U'\colon L^p(\widetilde{\Omega})\to L^p(\widetilde{\Omega})$
be defined by $U'(z,w)=(Uz,w)$ and let $Q'\colon L^p(\widetilde{\Omega})\to L^p(\Omega)$
be defined by $Q'(z,w)=Q(z)$. Then $U'$ is an isomorphism, $({U'}^n)_{n\in\Zdb}$ is bounded
and $Q'{U'}^nJ'=T^n$ for any $n\geq 0$.
\end{proof}

\medskip
\section{Comparing $p$-boundedness properties}
In this section we will consider an $L^p$-analog of complete polynomial boundedness
going back to \cite{P2} (see also \cite[Chap. 8]{P4}) 
and give complements to the results obtained
in the previous section. In particular we will show the existence 
of $p$-polynomially bounded operators $L^p\to L^p$ without 
any loose dilation.

In the sequel we assume that $1\leq p<\infty$.
Let $n\geq 1$ be an integer. For any vector space $V$, we let $M_n(V)$
denote the space of $n\times n$ matrices with entries in $V$. When
$V=B(X)$ for some Banach space $X$, we equip this space with a specific
norm, as follows. For any $[T_{ij}]_{1\leq i,j\leq n}$ in
$M_n\bigl(B(X)\bigr)$, we set
\begin{equation}\label{5Pnorm}
\bignorm{[T_{ij}]}_{p,M_n(B(X))}\, =\, \sup\Biggl\{\biggl(
\sum_{i=1}^{n} \biggnorm{\sum_{j=1}^{n}
T_{ij}(x_j)}_{X}^p\biggr)^{\frac{1}{p}}\, : \, x_1,\ldots, x_n\in
X,\ \sum_{j=1}^{n} \norm{x_j}^p_{X} \leq 1 \Biggr\}.
\end{equation}
In other words, we regard $[T_{ij}]$ as an operator
$\ell^p_n(X)\to \ell^p_n(X)$ in
a natural way and the norm of the matrix is defined as the
corresponding operator norm.

Let $X,Y$ be two Banach spaces, let $V\subset B(X)$ be a subspace and
let $u\colon V\to B(Y)$ be a linear mapping. We say that $u$ is $p$-completely
bounded if there exists a constant $C\geq 0$ such that
$$
\bignorm{[u(T_{ij})]}_{p,M_n(B(Y))}\,\leq C\bignorm{[T_{ij}]}_{p,M_n(B(X))}
$$
for any $n\geq 1$ and any matrix $[T_{ij}]$ in $M_n(V)$. In this case,
we let $\norm{u}_{pcb}$ denote the smallest possible $C$.

Let us regard the vector space $\P$ of all complex polynomials as a subspace
of $B(\ell^p_\Zdb)$, by identifying any $\varphi\in \P$ with the operator
$\varphi(S)$. Accordingly for any
$[\varphi_{ij}]$ in $M_n(\P)$, we set
$$
\bignorm{[\varphi_{ij}]}_{p}\,=\,\bignorm{[\varphi_{ij}(S)]}_{p, M_n(B(\ell^p_{\mathbb{Z}}))}.
$$
This extends (\ref{4Pnorm}) to matrices. We say that a bounded operator
$T\colon Y\to Y$ is $p$-completely polynomially bounded if
the natural mapping $u\colon\P\to B(Y)$ given 
by $u(\varphi)=\varphi(T)$ is $p$-completely
bounded. This is equivalent to the existence
of a constant $C\geq 1$ such that
$$
\bignorm{[\varphi_{ij}(T)]}_{p,M_n(B(Y))}\, \leq C\bignorm{[\varphi_{ij}]}_{p}
$$
for any matrix $[\varphi_{ij}]$ of complex polynomials.

When $p=2$ and $Y$ is a Hilbert space, the notions of $2$-polynomial boundedness and $2$-complete polynomial boundedness
correspond to the usual notions of polynomial boundedness and complete
polynomial boundedness from \cite{Pa1,Pa2}. 
See \cite{Pa2}
for the rich connections with operator space theory. The existence of a polynomially
bounded operator on Hilbert space which is not completely
polynomially bounded is a major result due to Pisier. Indeed this
is the heart of his negative solution to
the Halmos problem \cite{P3,P4}. We will show that Pisier's
construction can be transferred to our $L^p$-setting.

We start with an elementary result which is obvious when $p=2$ but requires
attention when $p\not=2$.

\begin{lemma}\label{5MnRep}
Let $N\geq 1$ be an integer, let $H$ be a Hilbert space and
let $\pi\colon B(\ell^2_N)\to B(H)$ be a unital $*$-representation.
Then for any $n\geq 1$ and any matrix $[T_{ij}]$ in $M_n\bigl(B(\ell^2_N)\bigr)$, we have
$$
\bignorm{[T_{ij}]}_{p,M_n(B(\ell^2_N))}\,\leq\,
\bignorm{[\pi(T_{ij})]}_{p,M_n(B(H))}.
$$
\end{lemma}

\begin{proof} As is well-known, there is a Hilbert space
$K$ such that
$$
H\simeq \ell^2_N(K),\qquad B(H)\simeq B(\ell^2_N)\otimes B(K),
$$
and $\pi(T)=T\otimes I_K$
for any $T\in B(\ell^2_N)$ (see e.g. \cite[Cor. III.1.7]{Dav}). Consider  
$[T_{ij}]$ in $M_n\bigl(B(\ell^2_N)\bigr)$ and  $x_1,\ldots,x_n$ in
$\ell^2_N$. Fix some $e\in K$ with $\norm{e}=1$. Then 
\begin{align*}
\sum_i\Bignorm{\sum_{j} T_{ij}(x_j)}_{\ell^2_N}^p\, &
=\,\sum_i\Bignorm{\sum_{j} T_{ij}(x_j) \otimes e}_{\ell^2_N(K)}^p\\
&
=\,\sum_i\Bignorm{\sum_{j} \bigl[\pi(T_{ij})\bigr](x_j\otimes e)}_{\ell^2_N(K)}^p\\
&
\leq \,\bignorm{[\pi(T_{ij})]}^{p}_{p, M_n(B(H))} 
\,\sum_j\norm{x_j\otimes e}^{p}_{\ell^2_N(K)}\\
&
\leq \,
\bignorm{[\pi(T_{ij})]}^{p}_{p, M_n(B(H))}\,\sum_j\norm{x_j}^p_{\ell^2_N}\,,
\end{align*}
and the result follows at once.
\end{proof}

\begin{proposition}\label{5Pisier}
Suppose that $1<p<\infty$. There exists a $p$-polynomially bounded
operator $T\colon L^p\bigl([0,1]\bigr)\to L^p\bigl([0,1]\bigr)$ which is
not $p$-completely polynomially bounded.
\end{proposition}

\begin{proof} 
We need some background on Pisier's counterexample. We refer 
to \cite[Chap. 9]{P4} and \cite[Chap. 10]{Pa2} for a detailed exposition of this example 
and also to the necessary background on Hankel operators on $B(\ell^2(H))$
and their $B(H)$-valued symbols.

We start with a concrete description of a sequence of operators satisfying the
so-called canonical anticommutation relations.
Let $I_2$ denote the identity matrix on $M_2$. For any $k\geq 1$,
consider the unital embedding $M_{2^k}\hookrightarrow M_{2^{k+1}}\simeq 
M_{2^k}\otimes M_2$ given
by $A\mapsto A\otimes I_2$. The closure of the union of the resulting increasing sequence
$(M_{2^k})_{k\geq 1}$ is a $C^*$-algebra. Representing it as an algebra of operators, we
obtain a Hilbert space $H$ and an embedding
\begin{equation}\label{5Representation}
\mathop{\bigcup}\limits_{k\geq 1}\big\uparrow M_{2^k}\,\subset B(H)
\end{equation}
whose restriction to each $M_{2^k}$ is a unital $*$-representation.

Consider the $2\times 2$ matrices
$$
D=\left[
    \begin{array}{cc}
      0 & 1 \\
      0 & 0 \\
    \end{array}
  \right]
\qquad\text{and}\qquad 
E=\left[
    \begin{array}{cc}
      1 & 0 \\
      0 & -1 \\
    \end{array}
  \right].
$$
For any $k\geq 1$, we set 
$$
C_k=E^{\otimes (k-1)}\otimes D \ \in M_{2^k},
$$
where $E^{\otimes (k-1)}$ denotes the tensor product of $E$ with itself $(k-1)$ times.
Then following (\ref{5Representation}) 
we let $\widetilde{C_k}$ denote this operator regarded as an element of $B(H)$. The distinction between 
$C_k$ and $\widetilde{C_k}$ may look superfluous. The reason why we need this is that the
inclusion providing the identification between $C_k$ and $\widetilde{C_k}$
is a $*$-representation and a priori, $*$-representations are not $p$-complete isometries
(i.e. they do not preserve $p$-matrix norms). However using Lemma \ref{5MnRep}, we see that
for any $m\geq 1$, for any $n\geq 1$ and for any $a_1,\ldots,a_m\in M_{n}$,
\begin{equation}\label{5Increase}
\Bignorm{\sum_{k=1}^{m} a_k\otimes C_k\otimes I_2^{\otimes (m-k)}}_{p,M_n(B(\ell^2_{2^m}))}
\leq\Bignorm{\sum_{k=1}^{m} a_k\otimes \widetilde{C_k}}_{p, M_n(B(H))}.
\end{equation}

The above sequence of matrices has the following remarkable property (see \cite[p. 70]{P4}):
for any complex numbers
$\alpha_1,\ldots,\alpha_m$,
\begin{equation}
\label{5Clifford}
\bgnorm{\sum_{k=1}^{m}\alpha_k \,C_k\otimes I_2^{\otimes (m-k)}}_{B(\ell^{2}_{2^m})}=\bigg(\sum_{k=1}^{m}|\alpha_k|^2\bigg)^{\frac{1}{2}}.
\end{equation}

Let $\H=\ell^2(H)\mathop{\oplus}\limits^2\ell^2(H)$, let
$\sigma\colon\ell^2(H)\to \ell^2(H)$ denote the shift operator, let 
$\Gamma\colon \ell^2(H)\to \ell^2(H)$ be the Hankel operator associated
to the $B(H)$-valued function $F$ given by
$$
F(t)=\,\sum_{k=1}^{\infty} \frac{\widetilde{C_k}}{2^k} \, e^{-i(2^k-1)t},
$$
and let $T\in B(\H)$ be the operator given by 
$$
T=\left[
    \begin{array}{cc}
      \sigma^*   & \Gamma \\
            0        & \sigma \\
    \end{array}
  \right].
$$

Pisier proved that this operator is polynomially bounded without being completely polynomially
bounded. Since $\norm{\ }_2\leq\norm{\ }_p$ on $\P$, the linear mapping
$$
u\colon (\P,\norm{\ }_p)\longrightarrow B(\H),\qquad u(\varphi)=\varphi(T),
$$
is therefore bounded. Our aim is now to show that $u$ is not $p$-completely bounded.

We consider the auxiliary mapping $w\colon\P\to B(H)$ defined by letting
$$
w\biggl(\sum_{k\geq 0} d_k z^k\biggr)\,=\,\sum_{k\geq 0} d_{2^k}\widetilde{C_k}
$$
for any finite sequence $(d_k)_{k\geq 0}$ of complex numbers.

Let $j\colon H\to \ell^2(H)$ be the isometric embedding given by 
$j(x)=(x,0,\ldots,0,\ldots)$.
Then let $v\colon B(\ell^2(H))\to B(H)$ be defined by letting
$v(R)=j^*Rj$ for any $R\in B(H)$. It is easy to check that 
$v$ is $p$-completely bounded, with $\norm{v}_{pcb}=1$.
On the other hand, for any $\varphi\in \P$, we have
$$
\varphi(T)=\left[
\begin{array}{cc}
\varphi(\sigma^*)
& \Gamma\varphi'(\sigma) \\
 0         &   \varphi(\sigma)    \\
\end{array}
\right],
$$
see \cite[(9.7)]{P4}. Let $\widetilde{u}\colon\P\to B(\ell^2(H))$ be defined
by $\widetilde{u}(\varphi)=\Gamma\varphi'(\sigma)$. Then the argument in the proof of
\cite[Thm. 9.7]{P4}
shows that $w=v\widetilde{u}$. Thus if $u$ were $p$-completely bounded, 
then $w$ would be $p$-completely bounded as well. Let us show that this does 
not hold true.

Note that for any Banach space $X$, for any integer $N\geq 1$,
for any $T\in B(X)$ and for any $A\in B(\ell^1_N)$, we have
$$
\bignorm{A\otimes T\colon \ell^1_N(X)\longrightarrow \ell^1_N(X)}\,=\,\norm{A}_{B(\ell^1_N)}
\norm{T}_{B(X)}.
$$
This can be be seen as a consequence of the fact that $\ell^1_N(X)$ is 
the projective tensor product of $\ell^1_N$ and $X$, see \cite[Chapter VIII]{DU}, however 
an elementary proof is also possible (we leave this to the reader).

Let $m\geq 1$.
Clearly $\norm{E}_{B(\ell^1_2)}=\norm{D}_{B(\ell^1_2)}=1$. Hence applying the above
property we have $\norm{C_k}_{B(\ell^1_{2^k})}=\norm{E}_{B(\ell^1_2)}^{k-1}\norm{D}_{B(\ell^1_2)}=1$
and hence 
\begin{equation}\label{5p=1}
\bignorm{C_k\otimes I_2^{\otimes (m-k)}\otimes S^{2^k}}_{B(\ell^1_{2^m}(\ell^1_{\Zdb}))}\,
=1,\qquad k=1,\ldots,n.
\end{equation}
Let $\varphi_m\in M_{2^m}\otimes\P$ be given by
$$
\varphi_m(z)=\,\sum_{k=1}^{m} C_k\otimes I_2^{\otimes (m-k)}\, z^{2^k}.
$$
By (\ref{5Clifford}),  we have
$$
\norm{\varphi_m}_{2}=\sup_{|z|=1}\norm{\varphi_m(z)}_{B(\ell^2_{2^m})} =
\sup_{|z|=1}\bgnorm{\sum_{k=1}^m z^{2^k}\,C_k\otimes I_2^{\otimes (m-k)}}_{B(\ell^2_{2^m})} =\sqrt{m}.
$$
On the other hand, applying (\ref{5p=1}) we have
$$
\norm{\varphi_m}_1\,=\,
\biggnorm{\sum_{k=1}^m C_k\otimes I_2^{\otimes (m-k)}\otimes S^{2^k}}_{B(\ell^1_{2^m}(\ell^1_\Zdb))}
\,\leq\,\sum_{k=1}^m 
\bignorm{C_k\otimes I_2^{\otimes (m-k)}\otimes S^{2^k}}_{B(\ell^1_{2^m}(\ell^1_{\Zdb}))}\,
\,=\,m.
$$
By interpolation, we deduce that
\begin{equation}\label{5P-estimate}
\norm{\varphi_m}_{p} \leq m^{\frac{1}{p}}.
\end{equation}

Next we have
\begin{equation}\label{5w}
\bigl(I_{B(\ell^p_{2^m})}\otimes w\bigr) (\varphi_m)\,=\,\sum_{k=1}^m C_k\otimes I_2^{\otimes (m-k)}
\otimes \widetilde{C_k}\,.
\end{equation}
Let us estimate the norm of this tensor product in $B(\ell^p_{2^m}(H))$.  
Let $(e_1,e_2)$ denote the canonical basis of $\Cdb^2$. 
For any $k=1,\ldots, m$ and any
$i_1,\ldots,i_m,j_1,\ldots,j_m$ in $\{1,2\}$,
\begin{align*}
\MoveEqLeft\Bigl\langle \Bigl(C_k\otimes I_2^{\otimes (m-k)}\Bigr)
(e_{j_1}\ot\cdots\ot e_{j_m}),
e_{i_1}\ot\cdots\ot e_{i_m}\Bigr\rangle\\
    &=\bigl\langle E e_{j_1}\ot\cdots\ot Ee_{j_{k-1}}\ot De_{j_k}
    \ot e_{j_{k+1}}\cdots\ot e_{j_m},e_{i_1}\ot\cdots\ot e_{i_m}\bigr\rangle\\
    &=\bigl\langle (-1)^{\delta_{j_1,2}}e_{j_1}\ot\cdots\ot (-1)^{\delta_{j_{k-1},2}}e_{j_{k-1}}\ot \delta_{j_{k},2}e_{1}\ot e_{j_{k+1}}\cdots\ot e_{j_m},e_{i_1}\ot\cdots\ot e_{i_m}\bigr\rangle\\
    &=(-1)^{\delta_{j_1,2}}\cdots(-1)^{\delta_{j_{k-1},2}}
    \delta_{j_{k},2}\bigl\langle e_{j_1},e_{i_1}\bigr\rangle\cdots 
    \bigl\langle e_{j_{k-1}},e_{i_{k-1}}\bigr\rangle \bigl\langle e_{1},e_{i_{k}}\bigr\rangle\bigl\langle e_{j_{k+1}},e_{i_{k+1}}\bigr\rangle\cdots \bigl\langle e_{j_m},e_{i_m}\bigr\rangle\\ &=(-1)^{\delta_{j_1,2}}\cdots(-1)^{\delta_{j_{k-1},2}}\delta_{j_{k},2}
\delta_{1,i_{k}} \delta_{j_1,i_1}\cdots\delta_{j_{k-1},i_{k-1}}
\delta_{j_{k+1},i_{k+1}}\cdots\delta_{j_m,i_m}.
\end{align*}
Hence 
\begin{align*}
\biggl\langle\bigg(\sum_{k=1}^{m}\Bigl(C_k\otimes I_2^{\otimes (m-k)}\Bigr)
\otimes\Bigl(C_k\otimes I_2^{\otimes (m-k)}\Bigr)\biggr) & \biggl( \sum_{j_1,\ldots,j_m=1}^{2}
e_{j_1}\ot\cdots \ot e_{j_m}
\ot e_{j_1}\ot\cdots \ot e_{j_m}\biggr),\\
&
\sum_{i_1,\ldots,i_m=1}^{2} e_{i_1}\ot\cdots \ot e_{i_m}\ot e_{i_1}\ot\cdots \ot e_{i_m}\biggr\rangle
\end{align*}
is equal to
\begin{align*}
\ &\sum_{k=1}^m 
\sum_{\substack{j_1,\ldots,j_m,\\i_1,\ldots,i_m=1}}^{2}
\Bigl\langle \Bigl(C_k\otimes I_2^{\otimes (m-k)}\Bigr)
(e_{j_1}\ot\cdots\ot e_{j_m}),
e_{i_1}\ot\cdots\ot e_{i_m}\Bigr\rangle^2\\
= &
\sum_{k=1}^m 
\sum_{\substack{j_1,\ldots,j_m,\\i_1,\ldots,i_m=1}}^{2}
\bigl(\delta_{j_{k},2}
\delta_{1,i_{k}} \delta_{j_1,i_1}\cdots\delta_{j_{k-1},i_{k-1}}
\delta_{j_{k+1},i_{k+1}}\cdots\delta_{j_m,i_m}\bigr)^2\\
= &
\sum_{k=1}^m 2^{m-1}\, =\, m2^{m-1}.
\end{align*}
Since the norm of 
$$
\sum_{i_1,\ldots,i_m=1}^{2}e_{i_1}\ot\cdots \ot e_{i_m}\ot
e_{i_1}\ot\cdots \ot
e_{i_m}
$$
in $\ell^{p}_{2^m}(\ell^{2}_{2^m})$ (resp. in $\ell^{p'}_{2^m}(\ell^{2}_{2^m})$) is equal to 
$$
\biggl(\sum_{i_1,\ldots,i_m=1}^{2}
\norm{e_{i_1}\ot\cdots \ot
e_{i_m}}_{\ell^2_{2^m}}^p\biggr)^{\frac{1}{p}}\,
=\, 2^{\frac{m}{p}}
$$
(resp.  $2^{\frac{m}{p'}}$), we deduce that
$$
\biggnorm{\sum_{k=1}^{m}\Bigl(C_k\otimes I_2^{\otimes (m-k)}\Bigr)\otimes
\Bigl(C_k\otimes I_2^{\otimes (m-k)}\Bigr)}_{B(\ell^p_{2^m}(\ell^2_{2^m}))}
\geq \frac{m}{2}.
$$
Combining with (\ref{5Increase}) and (\ref{5w}) we obtain that
$$
\bignorm{\bigl(I_{B(\ell^{p}_{2^m})}\otimes w\bigr)(\varphi_m)}_{p,M_{2^m(B(H))}} 
\,\geq\,\frac{m}{2}\,.
$$
Together with (\ref{5P-estimate}), this implies that $w$ is not $p$-completely bounded.
Thus $T$ is not $p$-completely polynomially bounded.

So we are done except that $T$ acts on the Hilbert space $\H$ and not on
$L^p\bigr([0,1]\bigr)$. However arguing as in the last part of
the proof of Theorem \ref{4Dilatable},
it is easy to pass from $\H$ to the space
$L^p\bigr([0,1]\bigr)$.
\end{proof}

The proof of Proposition \ref{4Transference} actually yields
the following stronger result: if an operator $T\colon L^p(\Omega)\to L^p(\Omega)$ admits a
loose dilation, then it is $p$-completely polynomially bounded (details are left to the reader).
Hence the above proposition yields the following.

\begin{corollary}
There exists a $p$-polynomially bounded operator 
$T\colon L^p\bigl([0,1]\bigr)\to L^p\bigl([0,1]\bigr)$
which does not admit any loose dilation.
\end{corollary}

Note also that according to Theorem \ref{4Dilatable} and the above observation, no $R$-Ritt operator can satisfy 
Proposition \ref{5Pisier}. Namely, 
if $T\colon L^p(\Omega)\to L^p(\Omega)$ is an $R$-Ritt operator and is $p$-polynomially
bounded, then it is $p$-completely polynomially bounded.

Remark \ref{4p=2} and the above investigations lead to the
following open problem (for $p\not=2$): {\it does any $p$-completely polynomially bounded operator
$T\colon L^p(\Omega)\to L^p(\Omega)$ admit a loose dilation?}

\bigskip
In the last part of this section we are going to consider 
another type of counterexamples.
Clearly any $p$-polynomially bounded $T\colon L^p(\Omega)\to L^p(\Omega)$ 
is automatically power bounded, that is,
$$
\sup_{n\geq 0}\norm{T^n}\,<\infty\,.
$$
The existence of a power bounded operator on a Hilbert space which is
not polynomially bounded is an old result of Foguel \cite{Fo}. 
Our aim is to prove an $L^p$-analog of that result. 
We will actually show a stronger form: there exists a 
Ritt operator which is not 
$p$-polynomialy bounded. To
achieve this, we will 
adapt 
the approach used in
\cite{KL} to go beyond Foguel's Theorem. 

We need some background on Schauder bases and their
multipliers that we briefly recall.
We let 
$v_1$ denote the set of all sequences $(c_n)_{n\geq 0}$ of complex numbers
whose variation $\sum_{n=1}^\infty \vert c_{n}
-c_{n-1}\vert\,$ is finite. Any such sequence is bounded and
$v_1$ is a Banach space for the 
norm
$$
\bignorm{(c_n)_{n\geq 0}}_{v_1}\,=\,\vert c_0\vert\, +\,
\sum_{n=1}^\infty \vert c_{n}
-c_{n-1}\vert\,.
$$

Let $(e_n)_{n\geq 0}$ be a Schauder basis on some Banach space $X$.
For any $n\geq 0$, let $Q_n\colon X\to X$ be the projection defined
by
\begin{equation}\label{5Qn}
Q_n\biggl(\sum_{k=0}^{\infty} a_ke_k\biggr)\,=\,\sum_{k=n}^{\infty}  a_ke_k
\end{equation}
for any converging sequence $\sum_k a_k e_k$. The sequence $(Q_n)_{n\geq 0}$ is bounded
and by a standard Abel summation argument, we have the following.

\begin{lemma}\label{5Abel}
For any $c=(c_n)_{n\geq 0}$ in $v_1$, there exists a (necessarily unique)
bounded operator $T_c\colon X\to X$ such that
$$
T_c\biggl(\sum_{n=0}^{\infty} a_ne_n\biggr)\,=\,\sum_{n=0}^{\infty} c_n a_ne_n
$$
for any converging sequence $\sum_n a_n e_n$. Furthermore,
$$
\norm{T_c}\leq \Bigl(\sup_{n\geq 0}\norm{Q_n}\Bigr) \bignorm{(c_n)_{n\geq 1}}_{v_1}.
$$
\end{lemma}

The above operator $T_c$ is called the multiplier associated to the sequence $c$.

\begin{proposition}\label{5Foguel}
Let $1<p<\infty$.
\begin{itemize}
\item [(1)] There exists a Ritt (hence a power bounded) operator on $\ell^2$ 
which is not $p$-polynomially bounded.
\item [(2)] 
There exists an $R$-Ritt operator on  $L^p\bigl([0,1]\bigr)$ 
which is not $p$-polynomially bounded.
\end{itemize}
\end{proposition}

\begin{proof}
(1): 
We let $(e_n)_{n\geq 0}$ be a Schauder basis of $H=\ell^2$.
It is clear that the sequence $\bigl(1-\frac{1}{2^n}\bigr)_{n\geq 0}$
has a finite variation. According to the above discussion, we
let $T\colon H\to H$ denote the multiplier associated
to this sequence. 

For any $\theta\in (-\pi,0)\cup(0,\pi]$, set
$$
c(\theta)_{n}=\,\frac{1}{e^{i\theta} - \bigl(1-\frac{1}{2^n}\bigr)}\,,\qquad n\geq 0.
$$
We have 
\begin{align*}
\sum_{n=1}^{+\infty}\vert c(\theta)_{n} -c(\theta)_{n-1}\vert\,
&= \,\sum_{n=1}^{+\infty} \,\Bigl\vert \int_{1-\frac{1}{2^{n-1}}}^{1-\frac{1}{2^{n}}}
\frac{dt}{\bigl(e^{i\theta}-t\bigr)^{2}}\,\Bigr\vert \\
&\leq\, \sum_{n=1}^{+\infty} \int_{1-\frac{1}{2^{n-1}}}^{1-\frac{1}{2^{n}}} 
\frac{dt}{\vert e^{i\theta}-t\vert^{2}}\\
&\leq\, \int_{0}^{1} \frac{dt}{\vert e^{i\theta} -t\vert^{2}}\,.
\end{align*}
Let $I(\theta)$ denote the latter integral. It is finite hence 
$c(\theta)=(c(\theta)_{n})_{n\geq 0}$ belongs to $v_1$. It is easy to deduce
that $e^{i\theta} -T$ is invertible, the operator $R(e^{i\theta},T)$ being
the multiplier associated to the sequence $c(\theta)$. 

For $\theta\in (0,\pi)$, elementary
computations yield
\begin{align*}
I(\theta) & = \,\int_{0}^{1}\frac{dt}{(t-\cos(\theta))^2+\sin^2(\theta)}\\
& =\,\frac{1}{\sin(\theta)}\,
\int_{-\frac{\cos(\theta)}{\sin(\theta)}}^{\frac{1-\cos(\theta)}{\sin(\theta)}}
\frac{du}{1+u^2}\\
& =\, \frac{\pi-\theta}{2\sin(\theta)}\,.
\end{align*}
Moreover $\vert e^{i\theta}- 1 \vert =2\sin\bigl(\tfrac{\theta}{2}\bigr)$,
hence
$$
\vert e^{i\theta} -1\vert I(\theta) \, =\, (\pi-\theta)
\frac{\sin\bigl(\tfrac{\theta}{2}\bigr)}{\sin(\theta)}\, =
\,\frac{\pi-\theta}{2\cos\bigl(\tfrac{\theta}{2}\bigr)}\,.
$$
This is bounded for $\theta$ varying in $(0,\pi)$, and $I(\theta)=I(-\theta)$ when 
$\theta\in (-\pi,0)$. 
According to Lemma \ref{5Abel}, this shows that 
$$
\sigma(T)\subset\Ddb\cup\{1\}\qquad\hbox{and}\qquad
\bigl\{(\lambda-1) R(\lambda,T)\, :\,\lambda\in\Tdb\setminus\{1\}\bigr\}\quad\hbox{is bounded.}
$$
Applying the maximum principle to the function $z\mapsto (1-z)(I_H-zT)^{-1}$, we deduce that the set
$\{(\lambda-1)R(\lambda,T)\, :\,\vert\lambda\vert>1\}$ is bounded as well, and hence 
$T$ is a Ritt operator.

Let us now assume that the basis $(e_n)_{n\geq 0}$ is not
an unconditional one. The operator $I-T$ is the multiplier
associated to the sequence $\big(\frac{1}{2^n}\big)_{n\geq 0}$ and
as is well-known, the lack of unconditionality implies that
for any $\theta \in (0,\pi)$,
this operator does not have a bounded
$H^\infty(\Sigma_\theta)$ functional calculus 
(see e.g. \cite[Thm. 4.1]{LM4} and its proof). 
According to Proposition
\ref{4Fourier}, this implies that 
$T$ is not $p$-polynomially bounded. 

\smallskip 
(2): Since all bounded subsets of $B(\ell^2)$ are $R$-bounded, the operator
considered in part (1) is automatically an $R$-Ritt operator. Then arguing again
as in the proof of Theorem \ref{4Dilatable},
it is easy to pass from an $\ell^2$-operator to an
$L^p\bigr([0,1]\bigr)$-operator which is not $p$-polynomially bounded although
being an $R$-Ritt operator.
\end{proof}

\medskip
\section{Generalizations to other Banach spaces}
Up to now we have mostly dealt with operators acting on (commutative) $L^p$-spaces.
In this last section, we shall consider more general Banach spaces, in
particular noncommutative $L^p$-spaces. We aim at extending our main results
from Sections 3 and 4 to this broader context. 

We will use classical notions from Banach space theory such as cotype, $K$-convexity and the UMD property.
We refer the reader to \cite{Bu,Die,P0} for background.

In accordance with (\ref{Rad(Lp)}), we are going to extend the definitions (\ref{1SF}) and
(\ref{1Alpha}) to arbitrary Banach spaces using Rademacher averages. Recall Section
2 for notation. The use of such averages as a substitute of square 
functions on abstract Banach spaces is a classical and fruitful principle. 
See e.g. \cite{JLX, KKW, LM0}.

Let $X$ be a Banach space, let $T\colon X\to X$ be any bounded operator and
let $x\in X$. Consider the element $x_k=k^{\frac{1}{2}}(T^{k}(x) - T^{k-1}(x))$
for any $k\geq 1$. If the series $\sum_k \varepsilon_k\otimes x_k\,$ converges in $L^2(\Omega_0;X)$
then we set
$$
\norm{x}_{T,1}\, =\,\biggnorm{\sum_{k=1}^{\infty} k^{\frac{1}{2}} \,\varepsilon_k\otimes \bigl(T^{k}(x) - T^{k-1}(x)
\bigr)}_{\Rad(X)}.
$$
We set $\norm{x}_{T,1}=\infty$ otherwise. Likewise, if $T$ is a Ritt operator and $\alpha>0$
is a positive real number, then we set
$$
\norm{x}_{T,\alpha}\, =\,\biggnorm{\sum_{k=1}^{\infty}  k^{\alpha - 
\frac{1}{2}}\, \varepsilon_k\otimes T^{k-1}(I-T)^\alpha
x}_{\Rad(X)}
$$
if the corresponding series converges in $L^2(\Omega_0;X)$, and
$\norm{x}_{T,\alpha}=\infty$ otherwise. 
The following extends Theorem \ref{3Equiv}.

\begin{theorem}\label{6Equiv}
Assume that $X$ is reflexive and has a finite cotype. Let $T\colon X\to X$ be an $R$-Ritt
operator. Then for any $\alpha>0$ and $\beta>0$, we have an equivalence
$$
\norm{x}_{T,\alpha}\,\approx\norm{x}_{T,\beta},\qquad x\in X.
$$
\end{theorem}

\begin{proof}
We noticed in Section 2 that if $X$ has a finite cotype, then Rademacher
averages and Gaussian averages are equivalent on $X$. 

Furthermore, the reflexivity of $X$ ensures that it satisfies
the Mean Ergodic Theorem. We thus have
$$
X=\,{\rm Ker}(I-T)\oplus\overline{{\rm Ran}(I-T)}.
$$

Lastly, since $X$ has a finite cotype, it cannot contain $c_0$ (as an isomorphic subspace).
Hence by \cite{Kwa}, a series
$\sum_k \varepsilon_k\otimes x_k\,$ converges in $L^2(\Omega_0;X)$
if (and only if) its partial sums are uniformly bounded, that is, there is
a constant $K\geq 0$ such that 
$$
\biggnorm{\sum_{k=1}^{N} \varepsilon_k\otimes x_k}_{{\rm Rad}(X)}\,\leq K,
\qquad N\geq 1.
$$

With these three properties in hand, it is easy to see that our
proof of Theorem \ref{3Equiv} extends verbatim to the general case.
\end{proof}

In the rest of this section we are going to focus on noncommutative $L^p$-spaces.
We let $M$ be a semifinite von Neumann algebra equipped with a normal 
semifinite faithful trace and for any $1\leq p<\infty$, we let $L^p(M)$ denote the 
associated (noncommutative) $L^p$-space. We refer to \cite{PX} for background and information
on these spaces. Any element of $L^p(M)$ is a (possibly unbounded)
operator and for any such $x$, we set
$$
\vert x\vert =(x^*x)^{\frac{1}{2}}.
$$
We recall the noncommutative analog of (\ref{Rad(Lp)}) from \cite{LP} (see also \cite{LPP}). 
For finite families $(x_k)_k$ of $L^p(M)$, we have the following
equivalences. 
If $2\leq p<\infty$, then 
\begin{equation}\label{6K1}
\biggnorm{\sum_{k} \varepsilon_k\otimes x_k}_{\Rad(L^p(M))}\,
\approx \, \max\biggl\{\biggnorm{\biggl(\sum_{k}
\vert x_k\vert^2\biggr)^{\frac{1}{2}}}_{L^p(M)}\, ,\, 
\biggnorm{\biggl(\sum_{k}
\vert x^*_k\vert^2\biggr)^{\frac{1}{2}}}_{L^p(M)}\biggr\}.
\end{equation}
If $1<p\leq 2$, 
then 
\begin{equation}\label{6K2}
\biggnorm{\sum_{k} \varepsilon_k\otimes x_k}_{\Rad(L^p(M))}\,
\approx \, \inf\Biggl\{\biggnorm{\biggl(\sum_{k}
\vert u_k\vert^2\biggr)^{\frac{1}{2}}}_{L^p(M)}\,+\,
\biggnorm{\biggl(\sum_{k}
\vert v^*_k\vert^2\biggr)^{\frac{1}{2}}}_{L^p(M)}\Biggr\},
\end{equation}where the infimum runs over all possible decompositions
$x_k=u_k+v_k$ in $L^p(M)$.

Let $T\colon L^p(M)\to L^p(M)$ be a bounded operator. We say that $T$ admits
a noncommutative loose dilation if there exist a von Neumann algebra
$\widetilde{M}$, an isomorphism $U\colon L^p(\widetilde{M})\to L^p(\widetilde{M})$ such that 
the set $\{U^n\,:\, n\in\Zdb\}$ is bounded and 
two bounded maps $L^p(M)\mathop{\longrightarrow}
\limits^{J} L^p(\widetilde{M})$ and $L^p(\widetilde{M}) \mathop{\longrightarrow}\limits^{Q}L^p(M)$ 
such that $T^n=QU^nJ$ for any integer $n\geq 0$. We
say that $T$ admits
a noncommutative strict dilation if this holds true for an isometric
isomorphism $U$ and two contractions
$J$ and $Q$. As opposed to the commutative case (see Remark \ref{4Strict}), there is no characterization
of contractions $T\colon L^p(M)\to L^p(M)$ which admit a noncommutative strict dilation.
The gap with the commutative situation is illustrated by the following result \cite[Thm 5.1]{JLM}: 
for any $p\not =2$, there exist a completely positive contraction  on some finite
dimensional noncommutative $L^p$-space with no noncommutative strict dilation.

We now turn to loose dilations. In the commutative setting, the following
proposition is a combination
of Propositions \ref{4Transference} and \ref{4Fourier}.

\begin{proposition}\label{6UMD}
Let $T\colon L^p(M)\to L^p(M)$, with $1<p<\infty$. If 
$T$ admits a noncommutative loose dilation, then
$I-T$ is sectorial and admits a bounded $H^{\infty}(\Sigma_\theta)$ functional 
calculus for any $\theta\in 
\bigl(\frac{\pi}{2},\pi\bigr)$.
\end{proposition}

\begin{proof} Let us explain how to adapt the `commutative' proof to
the present setting.

First we extend the definition (\ref{4Pnorm}) as follows. 
For any Banach space $X$, let $S_X\colon \ell^p_\Zdb(X)\to\ell^p_\Zdb(X)$ denote the
shift operator. Then for any $\varphi\in\P$, we set
$$
\norm{\varphi}_{p,X}=\bignorm{\varphi(S_X)}_{B(\ell^p_\Zdb(X))}.
$$
It follows from \cite[Thm. 4.3]{Bl1} that if $X$ is UMD, then 
Theorem \ref{4Marci} holds as well for scalar valued Fourier multipliers on 
$\ell^p_\Zdb(X)$. In this case, the argument in the proof of Proposition \ref{4Fourier}
leads to the following: for any $\theta\in 
\bigl(\frac{\pi}{2},\pi\bigr)$, there is an estimate
\begin{equation}\label{6UMD1}
\norm{\varphi}_{p,X}\lesssim\norm{\varphi}_{H^{\infty}(\Ddb_\theta)}
\end{equation}
for rational functions $\varphi$ with poles outside $\overline{\Ddb_\theta}$.

Second we note that if $U\colon L^p(\widetilde{M})\to L^p(\widetilde{M})$ is an isomorphism
such that $K=\sup\{\norm{U^n}\, :\, n\in\Zdb\}<\infty$, then the vectorial version
of the transference principle (see \cite[Thm. 2.8]{BG0})
ensures that for any $\varphi$ as above, we have
$$
\norm{\varphi(U)}\leq K^2\norm{\varphi}_{p,L^p(\widetilde{M})}.
$$

Assume now that $T\colon L^p(M)\to L^p(M)$ admits a noncommutative loose dilation. 
Noncommutative $L^p$-spaces are UMD hence property (\ref{6UMD1}) applies to them.
Hence arguing as in the proof of Proposition \ref{4Transference}, we find an estimate
$$
\norm{\varphi(T)}\,\lesssim\,\norm{\varphi}_{H^{\infty}(\Ddb_\theta)}
$$
for rational functions $\varphi$ with poles outside $\overline{\Ddb_\theta}$.
Finally the argument at the end of the proof of Proposition \ref{4Fourier}
yields that $I-T$ admits a bounded $H^{\infty}(\Sigma_\theta)$ functional calculus 
for any $\theta> 
\frac{\pi}{2}$. We skip the details.
\end{proof}

We are now ready to give 
a noncommutative analog of Theorem \ref{4Dilatable}.

\begin{theorem}\label{6Dilatable}
Let $T\colon L^p(M)\to L^p(M)$ be 
an $R$-Ritt operator, with $1<p<\infty$.
\begin{enumerate}
\item [(1)] 
The following assertions are equivalent.
\begin{enumerate}
\item [(i)] The operator $T$ 
admits a noncommutative loose dilation.
\item [(ii)] The operator $T$ and its adjoint $T^*\colon L^{p'}(M)\to L^{p'}(M)$
both satisfy uniform estimates
$$
\norm{x}_{T,1}\,\lesssim\,\norm{x}_{L^p(M)}\qquad\hbox{and}\qquad
\norm{y}_{T^*,1}\,\lesssim\,\norm{y}_{L^{p'}(M)}
$$
for $x\in L^p(M)$ and $y\in L^{p'}(M)$.
\end{enumerate}
\item [(2)] Assume that $p\geq 2$. Then the above conditions are equivalent to
the existence of a constant $C\geq 1$ for which the following two properties hold.
\begin{enumerate}
\item [(iii)] For any $x\in L^p(M)$,
$$
\biggnorm{\Bigl( 
\sum_{k=1}^{\infty} k\bigl\vert 
T^k(x)-T^{k-1}(x)\bigr\vert^2\Bigr)^{\frac{1}{2}}}_{L^p(M)}\,\leq C\norm{x}_{L^p(M)}
$$
and 
$$
\biggnorm{\Bigl( 
\sum_{k=1}^{\infty} k\bigl\vert \bigl(T^k(x)-T^{k-1}(x)\bigr)^* 
\bigr\vert^2\Bigr)^{\frac{1}{2}}}_{L^p(M)}\,\leq C\norm{x}_{L^p(M)}.
$$
\item [(iii)$^*$] For any $y\in L^{p'}(M)$, there exist two sequences $(u_k)_{k\geq 1}$ and $(v_k)_{k\geq 1}$
of $L^{p'}(M)$ such that 
$$
\biggnorm{\Bigl(\sum_{k=1}^{\infty} \vert u_k\vert^2\Bigr)^{\frac{1}{2}}}_{L^{p'}(M)}\,
+\,
\biggnorm{\Bigl(\sum_{k=1}^{\infty} \vert v_k^*\vert^2\Bigr)^{\frac{1}{2}}}_{L^{p'}(M)}\,\leq C\norm{y}_{L^{p'}(M)}\,,
$$
and 
$$
u_k+v_k\,=\, k^{\frac{1}{2}}\bigl( T^{* k}(y)-T^{*(k-1)}(y)\bigr)\qquad\hbox{for any}\ k\geq 1.
$$
\end{enumerate}
\end{enumerate}
\end{theorem}

\begin{proof}
Theorem \ref{3CF} holds as well on noncommutative $L^p$-spaces 
for $R$-Ritt operators, by \cite{LM0}. 
Combining that result with Proposition \ref{6UMD}, we obtain that (i) implies (ii).

Assume (ii) and suppose for simplicity that $I-T$ is $1$-$1$ (the changes to treat the 
general case are minor ones). 
By Theorem \ref{6Equiv}, we have uniform estimates
$$
\norm{x}_{T,\frac{1}{2}}\,\lesssim\,\norm{x}_{L^p(M)}\qquad\hbox{and}\qquad
\norm{y}_{T^*,\frac{1}{2}}\,\lesssim\,\norm{y}_{L^{p'}(M)}
$$
for $x\in L^p(M)$ and $y\in L^{p'}(M)$. As in the proof of Theorem
\ref{4Dilatable}, we may therefore define 
$J_1\colon L^p(M) \to {\rm Rad}(L^p(M))$ and 
$J_2\colon L^{p'}(M) \to {\rm Rad}(L^{p'}(M))$ by
setting 
$$
J_1(x)=\,\sum_{k=1}^{\infty} \varepsilon_k\otimes T^{k-1}(I-T)^{\frac{1}{2}}x
\qquad\hbox{and}\qquad
J_2(y)=\,\sum_{k=1}^{\infty} \varepsilon_k\otimes T^{*(k-1)}(I-T^*)^{\frac{1}{2}}y
$$
for any $x\in L^{p}(M)$ and any $y\in L^{p'}(M)$. Since $L^p(M)$ is $K$-convex, we have 
a natural isomorphism
\begin{equation}\label{6K}
\bigr({\rm Rad}(L^p(M))\bigr)^*\, \approx\, {\rm Rad}(L^{p'}(M)).
\end{equation}
Hence one can consider the composition $J_2^*J_1$, it is equal to $(I+T)^{-1}$
and one obtains (i) by simply adapting the proof of Theorem
\ref{4Dilatable}.

Finally the equivalence between (ii) and (iii)+(iii)$^*$ follows from 
(\ref{6K1}) and (\ref{6K2}).
\end{proof}

Note that switching (iii) and (iii)$^*$, we find a version of (2) for the case $p\leq 2$.

\begin{remark}
\  

(1) Let $T\colon L^p(\Omega)\to L^p(\Omega)$ be an $R$-Ritt operator on some {\it commutative} $L^p$-space. Combining 
Theorems \ref{6Dilatable} and \ref{4Dilatable}, we find that $T$
admits a noncommutative loose dilation (if and) only if it 
admits a commutative one.

(2) Proposition \ref{3DbSFE} holds true on noncommutative $L^p$-spaces. 
The proof is similar, replacing the square function by the norm in ${\rm Rad}(L^p(M))$ and using (\ref{6K})
instead of the duality $L^{p}(\Omega;\ell^2)^* = L^{p'}(\Omega;\ell^2)$.
\end{remark}

We refer the reader to \cite{Ar0} for examples of operators
with a noncommutative strict dilation, and to the paper \cite{Ar} for more about square
functions associated to Ritt operators on noncommutative $L^p$-spaces.


\begin{thebibliography}{99}
\bibitem{Ak} M. Akcoglu, {\it A pointwise ergodic theorem in $L_p$-spaces}, 
Canad. J. Math.  27 (1975), no. 5, 1075-1082.
\bibitem{A} M. Akcoglu, and L. Sucheston, {\it Dilations of positive contractions
on $L_p$ spaces}, Canad.  Math. Bull. 20 (1977), 285-292.
\bibitem{Ar0} C. Arhancet, {\it On Matsaev's conjecture for c
ontractions on noncommutative $L_p$-spaces}, to appear in J. 
Operator Theory,  arXiv:1009.1292.
\bibitem{Ar} C. Arhancet, {\it Square functions for Ritt operators on noncommutative $L^p$-spaces},
to appear in Math. Scandinavica, arXiv:1107.3415.
\bibitem{BG0} E. Berkson, and T. A. Gillespie, {\it Generalized analyticity in UMD spaces}, Ark. Mat. 27 (1989), 1-14.
\bibitem{BG} E. Berkson, and T. A. Gillespie, {\it Spectral decompositions and
harmonic analysis on UMD Banach spaces}, Studia Math. 112 (1994), 13-49.
\bibitem{Bl1} S. Blunck, {\it Maximal regularity of discrete and continuous time evolution equations},
Studia Math. 146  (2001),  no. 2, 157-176.
\bibitem{Bl2} S. Blunck, {\it Analyticity and discrete maximal regularity on $L_p$-spaces},
J. Funct. Anal.  183  (2001),  211-230.
\bibitem{Bu} D. L. Burkholder, {\it Martingales and singular integrals in Banach spaces}, pp. 233-269 in 
``Handbook of the geometry of Banach spaces", Vol. I,  North-Holland, Amsterdam, 2001.
\bibitem{CPSW} P. Cl\'{e}ment, B. de Pagter, F. A. Sukochev, and H. Witvliet,
{\it Schauder decompositions and multiplier theorems}, Studia
Math. 138 (2000), 135-163.
\bibitem{CRW} R. Coifman, R. Rochberg, and G. Weiss, {\it Applications of transference: the
$L_p$ version of von Neumann's inequality and the Littlewood-Paley-Stein theory}, pp. 53-67
in ``Linear spaces and Approximation", Birkh\"auser, Basel, 1978.
\bibitem{CW} R. R. Coifman, and G. Weiss, {\it Transference methods in analysis},
CBMS 31, Amer. Math. Soc., 1977.
\bibitem{CDMY} M. Cowling, I. Doust, A. McIntosh, and A. Yagi, {\it Banach
space operators with a bounded $H^{\infty}$ functional calculus},
J. Aust. Math. Soc., Ser. A  60 (1996), 51-89.
\bibitem{Dav} K. Davidson, {\it $C^*$-algebras by example}, 
Fields Institute Monographs, 6. American Mathematical Society, Providence, RI, 1996. xiv+309 pp.
\bibitem{Die} J. Diestel, H. Jarchow, and A. Tonge, {\it Absolutely summing operators}, 
Cambridge Studies in Advanced Mathematics 43, Cambridge University Press, Cambridge, 1995. xvi+474 pp. 
\bibitem{DU} J. Diestel, and J. J. Uhl, {\it Vector measures}, Mathematical Surveys, 
No. 15, American Mathematical Society, Providence, R.I., 1977. xiii+322 pp.
\bibitem{Dr} S. W. Drury, {\it A counterexample to a conjecture of Matsaev}, Lin. Alg. and its Appl. 435 (2011), 323-329.
\bibitem{EG} R. E. Edwards, and G. I.  Gaudry, {\it Littlewood-Paley and multiplier theory},
Ergebnisse der Mathematik und ihrer Grenzgebiete, Band 90. Springer-Verlag, Berlin-New York, 1977. ix+212 pp.
\bibitem{Fo} S. R. Foguel, {\it A counterexample to a problem of Sz.-Nagy},
Proc. Amer. Math. Soc. 15 (1964), 788-790.
\bibitem{FW} A. M. Fr\"ohlich, and L. Weis, {\it $H^\infty$ calculus and dilations}, Bull. Soc. Math. France  134  (2006),  no. 4, 487-508.
\bibitem{Ha} M. Haase, {\it The functional calculus for sectorial operators}, 
Operator Theory: Advances and Applications, 169, Birkhäuser Verlag, Basel, 2006. xiv+392 pp.
\bibitem{JLM} M. Junge, and C. Le Merdy, {\it Dilations and rigid factorisations on noncommutative $L^p$-spaces}, 
J. Funct. Anal.  249  (2007),  no. 1, 220-252. 
\bibitem{JLX} M. Junge, C. Le Merdy, and Q. Xu, {\it $H^{\infty}$
functional calculus and square functions on noncommutative
$L^p$-spaces}, Soc. Math. France, Ast\'erisque 305, 2006.
\bibitem{KKW} N. J. Kalton, Nigel, P. Kunstmann, and L. Weis, 
{\it Perturbation and interpolation theorems for the $H^\infty$-calculus with 
applications to differential operators},  
Math. Ann.  336  (2006), no. 4, 747-801.
\bibitem{KL} N. J. Kalton, and C. Le Merdy, 
{\it Solution of a problem of Peller concerning similarity},
J. Operator Theory 47 (2002), no. 2, 379-387.
\bibitem{KaW} N. J. Kalton, and L. Weis,
{\it The $H^\infty$ calculus and sums of closed operators}, Math. Annalen 321 (2001), 
319-345.
\bibitem{Kr} U. Krengel, {\it Ergodic theorems}, de Gruyter Studies in Mathematics, 6. Walter de Gruyter and
Co., Berlin, 1985. viii+357 pp.
\bibitem{KW} P. C. Kunstmann,  and L.  Weis, {\it Maximal $L_p$-regularity for
parabolic equations, Fourier multiplier theorems and
$H^\infty$-functional calculus}, pp. 65-311 in ``Functional
analytic methods for evolution equations", Lect. Notes in Math. 1855, Springer, 2004.
\bibitem{Kwa} S. Kwapien, {\it On Banach spaces containing $c_0$}, Studia Math. 52 (1974), 187-188.
\bibitem{LP} F. Lust-Piquard, {\it In\'egalit\'es de Khintchine dans $C_p\;(1<p<\infty)$} (French), 
C. R. Acad. Sci. Paris Sér. I Math. 303  (1986),  no. 7, 289-292. 
\bibitem{LPP} F. Lust-Piquard, and G. Pisier, {\it 
Noncommutative Khintchine and Paley inequalities},  Ark. Mat.  29  (1991),  no. 2, 241-260. 
\bibitem{LM0} C. Le Merdy, {\it $H^\infty$ functional calculus and square function estimates
for Ritt operators}, Preprint 2011,  arXiv:1202.0768. 
\bibitem{LM3} C. Le Merdy, {\it Square functions, bounded analytic semigroups, and applications},
Banach Center Publ.  75 (2007), 191-220.
\bibitem{LM4} C. Le Merdy, {\it $H^\infty$-functional calculus and applications to maximal regularity},  
Publ. Math. Besan\c con 16 (1998), 41-77.
\bibitem{LM1} C. Le Merdy,  {\it On square functions associated to sectorial operators},
Bull. Soc. Math. France  132  (2004),  no. 1, 137-156.
\bibitem{LMX1} C. Le Merdy, and Q. Xu, {\it Maximal theorems and square functions for analytic operators on 
$L^p$-spaces},  J. Lond. Math. Soc. (2) 86 (2012), no. 2, 343-365.
\bibitem{LMX2} C. Le Merdy, and Q. Xu, {\it Strong $q$-variation inequalities for analytic semigroups}, 
Annales Inst. Fourier 62 no. 6 (2012), 2069-2097.
\bibitem{Ly} Yu. Lyubich, {\it Spectral localization, power boundedness and invariant subspaces
under Ritt's type condition}, Studia Math. 134 (1999), 153-167.
\bibitem{MI} A. McIntosh, {\it Operators which have an $H^\infty$ functional calculus}, Proc. CMA Canberra
14 (1986), 210-231.
\bibitem{NZ} B. Nagy, and J. Zemanek, {\it A resolvent condition implying power boundedness}, Studia Math. 134 (1999), 143-151.
\bibitem{N} O. Nevanlinna, {\it Convergence of iterations for linear equations}, Birkha\"user, Basel, 1993.
\bibitem{Pa1} V. I. Paulsen, {\it Every completely polynomially bounded operator is similar to a contraction},
J. Funct. Anal.  55  (1984),  no. 1, 1-17.
\bibitem{Pa2} V. I. Paulsen, {\it Completely bounded maps and operator algebras},
Cambridge Studies in Advanced Mathematics, 78, Cambridge University Press, Cambridge, 2002. xii+300 pp. 
\bibitem{Pe} V. Peller, {\it An analogue of J. von Neumann's inequality for the space $L^{p}$} (Russian),
Dokl. Akad. Nauk SSSR  231  (1976), no. 3, 539-542.
\bibitem{P0} G. Pisier, {\it Probabilistic methods in the geometry of Banach spaces}, pp. 167–241 in
``Probability and analysis (Varenna, 1985)", Lecture Notes in Math. 1206, Springer, Berlin, 1986. 
\bibitem{P1} G. Pisier, {\it Complex interpolation and regular operators between Banach lattices}, 
Arch. Math. (Basel)  62  (1994),  no. 3, 261-269.
\bibitem{P2} G. Pisier, {\it Completely bounded maps between sets of Banach space operators},
Indiana Univ. Math. J.  39  (1990),  no. 1, 249-277.
\bibitem{P3} G. Pisier, {\it A polynomially bounded operator on Hilbert space which is not similar to a contraction},
J. Amer. Math. Soc.  10  (1997),  no. 2, 351-369.
\bibitem{P4} G. Pisier, {\it Similarity problems and completely bounded maps} (Second, expanded edition),
Lecture Notes in Mathematics, 1618. Springer-Verlag, Berlin, 1996. viii+156 pp.
\bibitem{P5} G. Pisier, {\it Complex Interpolation between Hilbert, Banach and Operator spaces},
Mem. Amer. Math. Soc.  208  (2010) vi+78 pp.
\bibitem{PX} G. Pisier, and Q. Xu, {\it Non-commutative
$L^p$-spaces}, pp. 1459-1517 in ``Handbook of the Geometry of
Banach Spaces", Vol. II, edited by W.B. Johnson and J. Lindenstrauss, Elsevier, 2003.
\bibitem{S2} E.M. Stein, {\it Topics in harmonic analysis related
to the Littlewood-Paley theory}, Ann. Math. Studies, Princeton,
University Press, 1970.
\bibitem{NF} B. Sz.-Nagy, and C. Foias, {\it Harmonic analysis of operators on Hilbert space}, Akademiai Kiad\'o,
Budapest, 1970.
\bibitem{V} P. Vitse, {\it A band limited and Besov class functional calculus for Tadmor-Ritt operators}, Archiv. Math.
 85  (2005),  no. 4, 374-385.
\end{thebibliography}
\end{document}